\documentclass[11pt, a4paper]{article}

\usepackage{amsmath}
\usepackage{amssymb}
\usepackage{amsthm}
\usepackage{amscd}
\usepackage{stmaryrd}
\usepackage{url}
\usepackage{tikz}
\usepackage{txfonts}
\usetikzlibrary{calc}

\usepackage[auth-sc]{authblk}

\usepackage[utf8]{inputenc}
\usepackage[english]{babel}
\usepackage[T1]{fontenc}
\usepackage{lmodern}
\usepackage{microtype}
\usepackage[
    linktocpage=true,
    colorlinks=true,
    linkcolor=purple,
    citecolor=black,
    urlcolor=purple,
    bookmarksnumbered=true,
]{hyperref}

\usepackage[totalwidth=160mm, totalheight=247mm]{geometry}

\usepackage{caption}
\renewcommand{\leq}{\leqslant}
\renewcommand{\geq}{\geqslant}
\newcommand{\bbN}{\mathbb{N}}
\newcommand{\bbZ}{\mathbb{Z}}

\newcommand{\bbR}{\mathbb{R}}

\newcommand{\mcA}{\mathcal A}

\newcommand{\mcU}{\mathcal U}
\newcommand{\mcD}{\mathcal D}
\newcommand{\mcO}{\mathcal O}

\newcommand{\mcL}{\mathcal L}
\newcommand{\LAR}{\mathcal L_\textup{AR}}
\newcommand{\Ledge}{\mathcal L_\textup{edge}}
\newcommand{\omcL}{\overline{\mathcal L}}
\newcommand{\oLAR}{\overline{\mathcal L_\textup{AR}}}
\newcommand{\oLedge}{\overline{\mathcal L_\textup{edge}}}

\newcommand{\bfe}{{\mathbf e}}
\newcommand{\bfv}{{\mathbf v}}
\newcommand{\bfx}{{\mathbf x}}
\newcommand{\bfy}{{\mathbf y}}
\newcommand{\Gv}{{\Gamma_{\bf v}}}
\newcommand{\Gvb}{{\Gamma_\vb}}
\newcommand{\transp}[1]{{{}^\textup{t} #1}}
\newcommand{\EOS}{{{\mathbf E}_1^*}}
\newcommand{\EOSS}{{{\mathbf E}_1^*(\sigma)}}

\newcommand{\ub}{{{\mathbf u}_\beta}}
\newcommand{\vb}{{{\mathbf v}_\beta}}
\newcommand{\Ms}{{{\mathbf M}_\sigma}}
\newcommand{\Msinv}{{{\mathbf M}^{-1}_\sigma}}
\newcommand{\Pc}{{\mathcal P_{\textup c}}}

\newcommand{\pic}{{\pi_{\textup{c}}}}

\newcommand{\myvcenter}[1]{\ensuremath{\vcenter{\hbox{#1}}}}

\newcommand{\svect}[3]{%
\big(\begin{smallmatrix}%
#1 \\ #2 \\ #3%
\end{smallmatrix}\big)%
}

\theoremstyle{plain}
\newtheorem{theo}{Theorem}[section]
\newtheorem{prop}[theo]{Proposition}
\newtheorem{lemm}[theo]{Lemma}
\newtheorem{coro}[theo]{Corollary}
\theoremstyle{definition}
\newtheorem{defi}[theo]{Definition}
\newtheorem{exam}[theo]{Example}
\newtheorem{rema}[theo]{Remark}
\newtheorem*{ackn}{Acknowledgements}

\title{\textbf{Substitutive Arnoux-Rauzy sequences \\ have pure discrete spectrum}}

\author[1]{Val\'erie Berth\'e}
\author[1,2]{Timo Jolivet}
\author[3]{Anne Siegel}

\affil[1]{
    LIAFA,
    Universit\'e Paris Diderot,
    France
}
\affil[2]{
    FUNDIM,
    Department of Mathematics,
    University of Turku,
    Finland
}
\affil[3]{
    IRISA,
    Campus de Beaulieu, Rennes,
    France
}

\date{}

\begin{document}
\maketitle

\begin{abstract}
We prove that the symbolic dynamical system generated by a purely substitutive Arnoux-Rauzy sequence
is measurably conjugate to a toral translation.
The proof is based on an explicit construction of a fundamental domain with fractal boundary
(a Rauzy fractal) for this toral translation.
\end{abstract}

\begin{center}
\emph{Dedicated to the memory of G\'erard Rauzy.}
\end{center}

\section{Introduction}
\paragraph{Substitutions and self-similarity}
Self-similarity is a notion that has proved its relevance in the study of dynamical systems.
Among the numerous ways of producing self-similar dynamical systems,
substitutions play a prominent role.
The term ``substitution'' encompasses
symbolic replacement rules that act on words (generating symbolic dynamical systems),
and more generally replacement rules acting on finite sets of compact subsets of $\bbR^d$
(generating tiling spaces, the most celebrated example being the Penrose tiling).
Substitutions prove to be very useful in many mathematical fields
(combinatorics on words, number theory, harmonic analysis, ergodic theory),
as well as in theoretical computer science and physics,
such as illustrated by the bibliography of \cite{PF02}.
 
\paragraph{$S$-adic systems}
A natural generalization of substitutive symbolic dynamical systems
arises when one considers the dynamical system obtained by iterating not one,
but several substitutions taken from a finite set $S$ of substitutions:
such systems are said \emph{$S$-adic}.
For more details see, \emph{e.g.}, \cite{Dur00, Leroy11}.
These systems, which have hierarchical structure,
belong to the larger category of fusion tilings introduced  in \cite{PS11},
which also includes Bratteli-Vershik systems as well as multi-dimensional cut-and-stack transformations. 
These hierarchies are not necessarily the same at each level
(as it is the case for purely substitutive systems),
but there are only finitely many different possible structures at each level.
The notion of $S$-adic system has been used notably for the description of linearly recurrent systems. 
For instance, Durand characterized linearly recurrent subshifts
as the set of primitive proper $S$-adic subshifts \cite{Dur00, Dur03}.

\paragraph{Arnoux-Rauzy sequences}
Sturmian sequences are among the most classical examples of $S$-adic sequences
and are a classical object of combinatorics on words and symbolic dynamics.
They are the infinite sequences of two letters with factor complexity $n+1$
(\emph{i.e.}, they have exactly $n+1$ factors of length $n$),
and they correspond to natural codings of irrational rotations on the circle~\cite{MH40}.
These sequences admit an $S$-adic definition:
they are the sequences that belong to the symbolic dynamical systems generated by iterating the two substitutions
$\sigma_1 : 1 \mapsto 1, 2 \mapsto 21$, and $\sigma_2 : 1 \mapsto 12, 2 \mapsto 2$.
See~\cite[Chap. 6]{PF02} and~\cite[Chap. 2]{Lot97} for a detailed survey of their properties.

Arnoux-Rauzy sequences were introduced in~\cite{AR91} in order to generalize
Sturmian sequences to three-letter alphabets.
They can be defined in an $S$-adic way as
the sequences that belong to the symbolic dynamical system generated by iterating the \emph{Arnoux-Rauzy substitutions}
$\sigma_1$, $\sigma_2$, $\sigma_3$ given by
\[
\sigma_i \ : \ \left\{
 \begin{array}{rcll}
 j &\mapsto& j & \text{ if } \ j = i \\
 j &\mapsto& ji & \text{ if } \ j \neq i
 \end{array}
\right.
\quad (i = 1,2,3),
\]
where each of the $\sigma_i$ occurs infinitely often in the iteration.
These sequences have factor complexity $2n+1$.
It was conjectured that Arnoux-Rauzy sequences correspond to natural codings
of translations on the two-dimensional torus
(as in the Sturmian case with rotations on the circle),
but this conjecture has been disproved in~\cite{CFZ00}.
This family of sequences has indeed a complex and intriguing ergodic behavior,
as underlined in \cite{CFM08}: some of them might even be measure-theoretically weak mixing. 
Arnoux-Rauzy sequences are also studied in word combinatorics,
they belong indeed to the family of episturmian words (see the survey \cite{AJ09}). 
For more references about Arnoux-Rauzy sequences, see~\cite{BFZ05, CC06}.
 
In this article, we focus on \emph{purely substitutive} Arnoux-Rauzy systems, that is,
systems generated by an infinite purely periodic product of Arnoux-Rauzy substitutions
(there are infinitely many such sequences).
We focus on the spectrum of these dynamical systems, and our main theorem is the following
(see Section \ref{sec:def} for precise definitions).

\begin{theo}
\label{theo:intro}
Let $\sigma$ be a finite product of the substitutions $\sigma_1$, $\sigma_2$, $\sigma_3$, 
where each substitution appears at least once.
The symbolic dynamical system $(X_\sigma, S)$ generated by $\sigma$ 
has pure discrete spectrum.
\end{theo}

\paragraph{Our methods}
In order to prove Theorem \ref{theo:intro},
we provide an explicit measure-theoretical isomorphism between
the symbolic dynamical system generated by an Arnoux-Rauzy substitution $\sigma$
and a two-dimensional translation of $\mathbb{T}^{2}$,
with the fundamental domain of the translation being the \emph{Rauzy fractal} of the substitution $\sigma$.
(See Section~\ref{subsec:RF} for more details about Rauzy fractals). 

In this article, we exploit a combinatorial definition of Rauzy fractals
introduced by Arnoux and Ito in \cite{AI01} to obtain the discrete spectrum property.
More precisely, we use methods inspired from \cite{IO94}:
we prove that a particular sequence of approximations of the Rauzy fractal contains arbitrarily large balls,
which is equivalent to the pure discrete spectrum property (see Section \ref{subsec:dynconseq}).

Our approach is well-suited for handling \emph{infinite} families of substitutions (as in Theorem \ref{theo:intro})
instead of a single given substitution.
Indeed, determining if a Pisot irreducible substitution has pure discrete spectrum is decidable \cite{Sie04, ST10, SS02,AL11},
but such decidability results are much more difficult to obtain for infinite families of substitutions,
although some results have been obtained in this direction in \cite{Akiya:2000,BBK06,BSW11}.
Our approach can also be applied to $S$-adic systems (see Theorem \ref{theo:growingballsSadic}).
We thus plan to develop the methods described here to further $S$-adic families of substitutions,
such as Jacobi-Perron substitutions in \cite{BBJS11}.

Furthermore, the present approach provides geometric information on the associated Rauzy fractals (see Theorem \ref{theo:RFprop} below and \cite{BJS11}).
This has interesting consequences for exhibiting Markov partitions of toral automorphisms.
Indeed, by Theorem \ref{theo:intro},
Rauzy fractals can be used to provide
explicit Markov partitions for the automorphism of the three-dimensional torus
defined by the abelianization of finite products of Arnoux-Rauzy substitutions
\cite{IO93, Pra99, IR06}.
Such Markov partitions were already known to exist \cite{Bow75},
but no general explicit construction is known in dimension $\geq 3$.
Note that the elements of the Markov partitions associated with Arnoux-Rauzy substitutions are connected,
by connectedness of their Rauzy fractals (see \cite{BJS11}).

Theorem \ref{theo:intro} has been stated in \cite{AI01} without proof,
and has recently been obtained independently by Barge, \v{S}timac and Williams \cite{BSW11}
using different methods based on the balanced pair algorithm and overlap coincidences.

\paragraph{Outline of the paper}
Section \ref{sec:def} introduces basic notions and definitions.
Specific combinatorial tools that allow the study of the topology
of the patterns generated by Arnoux-Rauzy dual substitutions
are described in Section \ref{sec:cov},
and in Section \ref{sec:mainresults},
we exploit the tools developed in the previous sections to prove our main results concerning the associated Rauzy fractals.

\section{Definitions and preliminary results}
\label{sec:def}
In this section, we introduce the basic objects of this article.
We start with the classical notion of substitutions on words,
and then we define unit faces and discrete planes. 
This allows us to define \emph{dual} substitutions,
which act on discrete planes instead of unidimensional words.
We then introduce Arnoux-Rauzy substitutions and Rauzy fractals,
together with basic considerations on the spectrum of symbolic dynamical systems.

\subsection{Substitutions}
Let $\mcA = \{1, \ldots, d\}$ be a finite alphabet,
and $\mcA^*$ the free monoid generated by $ \mcA$,
that is, the set of finite words over $\mcA$.
\begin{defi}[Substitution]
\label{defi:sub}
A \emph{substitution} is a non-erasing morphism of the free monoid $\mcA^*$,
\emph{i.e.}, a function $\sigma : \mcA^* \rightarrow \mcA^*$
such that $\sigma(uv) = \sigma(u)\sigma(v)$ for all words $u, v \in \mcA^*$,
and such that the image of each letter of $\mcA$ is non-empty.
\end{defi}
The \emph{incidence matrix} $\Ms$ of $\sigma$
is the square matrix of size $d \times d$ defined by $\Ms = (m_{ij})$,
where $m_{ij}$ is the number of occurrences of the letter $i$ in $\sigma(j)$.
A substitution $\sigma$ is \emph{unimodular} if $\det(\Ms) = \pm 1$.

A substitution $\sigma$ is \emph{primitive}
if there exists an integer $n$ (independent of the letters) such that 
$\sigma^n(a)$ contains at least one occurrence of the letter
$b$ for all $a,b \in \mcA$.
We will work with primitive substitutions having furthermore additional arithmetic properties.
We recall that a \emph{Pisot number} is a real algebraic integer
greater than $1$ whose conjugates have absolute value less than $1$.

\begin{defi}[Pisot irreducible substitution]
\label{def:Pisot}
A unimodular substitution $\sigma$ is \emph{Pisot irreducible}
if the maximal eigenvalue of $\Ms$ is a Pisot number,
and if the characteristic polynomial of $\Ms$ is irreducible.
\end{defi}

Let $\sigma$ be a unimodular Pisot irreducible substitution.
We denote by $\ub$ a positive eigenvector of $\Ms$ associated with $\beta$,
and by $\vb$ a positive eigenvector of $\transp \Ms$ associated with $\beta$.
The existence of such vectors comes from the Perron-Frobenius theorem applied to $\Ms$,
because the matrix of a Pisot irreducible substitution
is always primitive, as proved in~\cite{CS01}.

\begin{defi}[Contracting plane]
\label{defi:contrac}
Let $\sigma$ be a unimodular Pisot irreducible substitution.
The \emph{contracting plane} $\Pc$ associated with $\sigma$
is the plane of normal vector $\vb$.
We denote by $\pic : \bbR^3 \rightarrow \Pc$ the projection
of $\bbR^3$ on $\Pc$ along $\ub$.
\end{defi}

\paragraph{Symbolic dynamics of substitutions}
Let $S$ stand for the (two-sided) \emph{shift} on $\mcA^\bbZ$,
\emph{i.e.}, $S((u_n)_{n \in \bbZ})=(u_{n+1})_{n \in \bbZ}$.
We endow $\mcA^\bbZ$ with the product topology of the discrete topology on each copy of $\mcA$.
Let us now see how to associate with a primitive substitution a symbolic dynamical system,
defined as a closed shift invariant subset of $\mcA^\bbZ$.
Let $\sigma$ be a primitive substitution over ${\mathcal A}$,
and let $u \in \mcA ^\bbZ$ be such that $\sigma^k(u) = u$ for some $k \geq 1$
(such an infinite word exists by primitivity of $\sigma$).
Let $\overline{\mcO(u)}$ be the positive orbit closure of $u$ under the action of $S$,
\emph{i.e.}, the closure in $\mcA ^\bbZ$ of the set $\mcO(u) = \{S^n(u) : n \geq 0\}$.
By primitivity of $\sigma$, $\overline{\mcO(u)}$ does not depend on $u$,
so we can define $\overline{\mcO(u)}$ to be the \emph{symbolic dynamical system generated by $\sigma$},
that we denote by $X_\sigma$.
The system $(X_\sigma,S)$ is \emph{minimal} (every non-empty closed shift-invariant subset equals the whole set)
and \emph{uniquely ergodic} (there exists a unique shift-invariant probability measure on $X_\sigma$).
It consists of all the bi-infinite words
that have the same language as $u$.
See \cite{Que10} for more details.

\subsection{Discrete planes}
Let $(\bfe_1, \bfe_2, \bfe_3)$ denote the canonical basis of $\bbR^3$,
and let $\bbR_{>0}^3$ denote the set of vectors of $\bbR^3$ with positive entries.
The following definition of a discrete plane was introduced in \cite{IO93, IO94} under the name of ``stepped surface''.
The set of its vertices corresponds to the classic notion of ``standard arithmetic discrete plane'' in discrete geometry,
as introduced in \cite{Rev91, And03}.

\begin{defi}[Discrete plane]
\label{defi:discreteplane}
Let $\bfv \in \bbR_{>0}^3$ and
let $\mathcal S$ be the union of the unit cubes with integer coordinates
that are included in the upper half-space $\{\bfx \in \bbR^3 : \langle \bfx, \bfv \rangle \geq 0\}$.
The \emph{discrete plane $\Gv$ of normal vector $\bfv$} is the boundary of $\mathcal S$.
\end{defi}

A discrete plane can be seen as a union of \emph{unit faces} $[\bfx, i]^*$, defined by
\[
\definecolor{facecolor}{rgb}{0.8,0.8,0.8}
\begin{array}{ccccc}
  \,[\bfx, 1]^* & = & \{\bfx + \lambda \bfe_2 + \mu \bfe_3 : \lambda, \mu \in [0,1]\} & = &
    \myvcenter{%
    \begin{tikzpicture}
    [x={(-0.216506cm,-0.125000cm)}, y={(0.216506cm,-0.125000cm)}, z={(0.000000cm,0.250000cm)}]
    \fill[fill=facecolor, draw=black, shift={(0,0,0)}]
    (0, 0, 0) -- (0, 1, 0) -- (0, 1, 1) -- (0, 0, 1) -- cycle;
    \node[circle,fill=black,draw=black,minimum size=1.2mm,inner sep=0pt] at (0,0,0) {};
    \end{tikzpicture}} \\
  \,[\bfx, 2]^* & = & \{\bfx + \lambda \bfe_1 + \mu \bfe_3 : \lambda, \mu \in [0,1]\} & = &
    \myvcenter{%
    \begin{tikzpicture}
    [x={(-0.216506cm,-0.125000cm)}, y={(0.216506cm,-0.125000cm)}, z={(0.000000cm,0.250000cm)}]
    \fill[fill=facecolor, draw=black, shift={(0,0,0)}]
    (0, 0, 0) -- (0, 0, 1) -- (1, 0, 1) -- (1, 0, 0) -- cycle;
    \node[circle,fill=black,draw=black,minimum size=1.2mm,inner sep=0pt] at (0,0,0) {};
    \end{tikzpicture}} \\
  \,[\bfx, 3]^* & = & \{\bfx + \lambda \bfe_1 + \mu \bfe_2 : \lambda, \mu \in [0,1]\} & = & 
    \myvcenter{%
    \begin{tikzpicture}
    [x={(-0.216506cm,-0.125000cm)}, y={(0.216506cm,-0.125000cm)}, z={(0.000000cm,0.250000cm)}]
    \fill[fill=facecolor, draw=black, shift={(0,0,0)}]
    (0, 0, 0) -- (1, 0, 0) -- (1, 1, 0) -- (0, 1, 0) -- cycle;
    \node[circle,fill=black,draw=black,minimum size=1.2mm,inner sep=0pt] at (0,0,0) {};
    \end{tikzpicture}}%
\end{array}
\]
where $i \in \{1,2,3\}$ is the \emph{type}
and $\bfx \in \bbZ^3$ is the \emph{position} of $[\bfx, i]^*$, as illustrated:
\begin{center}
\input{fig/faces.tex}
\end{center}
We denote by $\bfx + f$ the translation of a face $f$ by $\bfx$,
and we denote by $\bfx + P$ the set $\{\bfx + f : f \in P\}$,
where $P$ is a union of faces and $\bfx \in \bbZ^3$.

The following proposition is a useful characterization of discrete planes,
that will often be used to prove that a given finite union of faces is not included in a given discrete plane.

\begin{prop}[\cite{ABI02, ABS04}]
\label{prop:arithplane}
The discrete plane $\Gv$ is the union of faces $[\bfx, i]^*$ satisfying
$0 \leq \langle \bfx, \bfv \rangle < \langle \bfe_i, \bfv \rangle$.
\end{prop}

\begin{rema}
\label{rema:imgU}
Every discrete plane contains the ``lower unit cube'' defined by
$\mcU = [\mathbf{0},1]^* \cup [\mathbf{0},2]^* \cup [\mathbf{0},3]^* =
\myvcenter{\begin{tikzpicture}
[x={(-0.216506cm,-0.125000cm)}, y={(0.216506cm,-0.125000cm)}, z={(0.000000cm,0.250000cm)}]
\definecolor{facecolor}{rgb}{0.8,0.8,0.8}
\fill[fill=facecolor, draw=black, shift={(0,0,0)}]
(0, 0, 0) -- (0, 0, 1) -- (1, 0, 1) -- (1, 0, 0) -- cycle;
\fill[fill=facecolor, draw=black, shift={(0,0,0)}]
(0, 0, 0) -- (0, 1, 0) -- (0, 1, 1) -- (0, 0, 1) -- cycle;
\fill[fill=facecolor, draw=black, shift={(0,0,0)}]
(0, 0, 0) -- (1, 0, 0) -- (1, 1, 0) -- (0, 1, 0) -- cycle;
\node[circle,fill=black,draw=black,minimum size=1.2mm,inner sep=0pt] at (0,0,0) {};
\end{tikzpicture}}$\,.
\end{rema}

\subsection{Dual substitutions}
We now introduce \emph{dual substitutions},
which act on unit faces of $\bbR^3$ instead of letters.
This notion, which can be seen as a multidimensional counterpart of substitutions,
was introduced by Ito and Ohtsuki~\cite{IO93, IO94} and refined later by Arnoux and Ito~\cite{AI01}.
A more general definition is given in~\cite{SAI01}.

\begin{defi}[Dual substitution]
\label{defi:dualsub}
Let $\sigma : \{1,2,3\}^* \rightarrow \{1,2,3\}^*$ be a unimodular substitution.
The \emph{dual substitution} $\EOSS$ associated with $\sigma$ is defined~by
\[
\EOSS([\bfx, i]^*) \ = \
\Msinv\bfx + 
\bigcup_{j = 1,2,3}
\;
\bigcup_{s \in \mcA^* \, : \, \sigma(j) = pis}
[\Msinv\ell(s), j]^*,
\]
where $\ell(s) = (|s|_1,|s|_2, |s|_3) \in \bbZ^3$ is the \emph{abelianization} of $s$
($|s|_i$ is the number of occurrences of the letter $i$ in $s$).
The notation ``$s \in \mcA : \sigma(j) = pis$'' in the above formula means that
the second union is done over all the words $s \in \mcA^*$ such that
$\sigma(j) = pis$ (with $p \in \mcA^*$),
\emph{i.e.,} over all the suffixes $s$ of $\sigma(j)$.
We extend the definition of $\EOSS$ to any union of unit faces:
$\EOSS(P_1 \cup P_2) \ = \ \EOSS(P_1) \cup \EOSS(P_2)$.
See Section \ref{subsec:AR} for some examples of dual substitutions.
\end{defi}

Let us remark that $\EOSS$ is completely described by $\Ms$ and
by the images of the faces $[\mathbf{0}, 1]^*$, $[\mathbf{0}, 2]^*$ and $[\mathbf{0}, 3]^*$,
by linearity of the multiplication by $\Msinv$.
Also, we have
$\EOS(\sigma \circ \sigma') = \EOS(\sigma') \circ \EOS(\sigma)$
for every unimodular substitutions $\sigma$ and $\sigma'$
(see~\cite{AI01} for more details).
The next proposition describes the interplay between discrete planes and dual substitutions.

\begin{prop}[\cite{AI01, Fer06}]
\label{prop:imgplane}
Let $\Gv$ be a discrete plane and $\sigma$ be a unimodular substitution.
The image of $\Gv$ by $\EOSS$ is the discrete plane $\Gamma_{\transp \Ms \bfv}$.

Furthermore,
two distinct faces of $\Gv$ have images by $\EOSS$
that do not have any unit face in common.
\end{prop}

Note that in the special case where $\Gv = \Gvb$
is the discrete plane associated with the contracting plane of $\sigma$
($\Gvb$ is indeed a discrete plane, because $\vb$ has positive coordinates),
Proposition \ref{prop:imgplane} implies that $\EOSS(\Gvb) = \Gamma_{\transp \Ms \vb} = \Gvb$,
which means that $\Gvb$ is invariant under $\EOSS$.

\subsection{Arnoux-Rauzy substitutions}
\label{subsec:AR}

Let $\sigma_1$, $\sigma_2$, $\sigma_3$ be the \emph{Arnoux-Rauzy substitutions} defined by
\[
\sigma_1 : \left\{
  \begin{array}{rcl}
  1 &\mapsto& 1 \\
  2 &\mapsto& 21 \\
  3 &\mapsto& 31
  \end{array}
\right.
\quad
\sigma_2 : \left\{
  \begin{array}{rcl}
  1 &\mapsto& 12 \\
  2 &\mapsto& 2 \\
  3 &\mapsto& 32
  \end{array}
\right.
\quad
\sigma_3 : \left\{
  \begin{array}{rcl}
  1 &\mapsto& 13 \\
  2 &\mapsto& 23 \\
  3 &\mapsto& 3
  \end{array}
\right..
\]
We are interested into the products of $\sigma_1, \sigma_2, \sigma_3$.
Such products verify some important properties
related with spectral properties addressed in Section \ref{subsec:RF}.

\begin{prop}[\cite{AI01}]
\label{prop:ARpisot}
Let $\sigma$ be a product of $\sigma_1$, $\sigma_2$ and $\sigma_3$
in which each substitution appears at least once.
Then, $\sigma$ is a Pisot irreducible substitution.
\end{prop}

Let $\Sigma_1$, $\Sigma_2$, $\Sigma_3$ be the dual substitutions defined by
$\Sigma_i = \EOS(\sigma_i)$ for $i \in \{1,2,3\}$.
According to Definition \ref{defi:dualsub}, we have
\[
\begin{array}{rcl}
\Sigma_1([\bfx, 1]^*) & = &
  \Msinv \bfx + [(0,0,0), 1]^* \cup [(0,0,0),2]^* \cup [(0,0,0), 3]^*\\
\Sigma_1([\bfx, 2]^*) & = &
  \Msinv\bfx + [(1,0,0), 2]^*\\
\Sigma_1([\bfx, 3]^*) & = &
  \Msinv\bfx + [(1,0,0), 3]^*
\end{array},
\]
\[
\begin{array}{rcl}
\Sigma_2([\bfx, 1]^*) & = &
  \Msinv \bfx + [(0,1,0), 1]^*\\
\Sigma_2([\bfx, 2]^*) & = &
  \Msinv\bfx +  [(0,0,0), 1]^* \cup [(0,0,0),2]^* \cup [(0,0,0), 3]^*\\
\Sigma_2([\bfx, 3]^*) & = &
  \Msinv\bfx + [(0,1,0), 3]^*
\end{array},
\]
\[
\begin{array}{rcl}
\Sigma_3([\bfx, 1]^*) & = &
  \Msinv \bfx +  [(0,0,1), 1]^*\\
\Sigma_3([\bfx, 2]^*) & = &
  \Msinv\bfx + [(0,0,1), 2]^*\\
\Sigma_3([\bfx, 3]^*) & = &
  \Msinv\bfx + [(0,0,0), 1]^* \cup [(0,0,0),2]^* \cup [(0,0,0), 3]^*
\end{array},
\]
which can be represented graphically as follows,
where the black dot stands for $\bfx$ on the left hand side
and for $\Msinv\bfx$ on the right-hand side.
\[
\definecolor{facecolor}{rgb}{0.8,0.8,0.8}
\renewcommand{\tabcolsep}{1.5mm}
\renewcommand{\arraystretch}{1.2}
\centering
\Sigma_1 : \left\{
    \begin{tabular}{rcl}%
    \myvcenter{%
    \begin{tikzpicture}
    [x={(-0.216506cm,-0.125000cm)}, y={(0.216506cm,-0.125000cm)}, z={(0.000000cm,0.250000cm)}]
    \fill[fill=facecolor, draw=black, shift={(0,0,0)}]
    (0, 0, 0) -- (0, 1, 0) -- (0, 1, 1) -- (0, 0, 1) -- cycle;
    \node[circle,fill=black,draw=black,minimum size=1.2mm,inner sep=0pt] at (0,0,0) {};
    \end{tikzpicture}}%
     & \myvcenter{$\mapsto$} & 
    \myvcenter{%
    \begin{tikzpicture}
    [x={(-0.216506cm,-0.125000cm)}, y={(0.216506cm,-0.125000cm)}, z={(0.000000cm,0.250000cm)}]
    \fill[fill=facecolor, draw=black, shift={(0,0,0)}]
    (0, 0, 0) -- (0, 1, 0) -- (0, 1, 1) -- (0, 0, 1) -- cycle;
    \fill[fill=facecolor, draw=black, shift={(0,0,0)}]
    (0, 0, 0) -- (0, 0, 1) -- (1, 0, 1) -- (1, 0, 0) -- cycle;
    \fill[fill=facecolor, draw=black, shift={(0,0,0)}]
    (0, 0, 0) -- (1, 0, 0) -- (1, 1, 0) -- (0, 1, 0) -- cycle;
    \node[circle,fill=black,draw=black,minimum size=1.2mm,inner sep=0pt] at (0,0,0) {};
    \end{tikzpicture}} \\
    \myvcenter{%
    \begin{tikzpicture}
    [x={(-0.216506cm,-0.125000cm)}, y={(0.216506cm,-0.125000cm)}, z={(0.000000cm,0.250000cm)}]
    \fill[fill=facecolor, draw=black, shift={(0,0,0)}]
    (0, 0, 0) -- (0, 0, 1) -- (1, 0, 1) -- (1, 0, 0) -- cycle;
    \node[circle,fill=black,draw=black,minimum size=1.2mm,inner sep=0pt] at (0,0,0) {};
    \end{tikzpicture}}%
     & \myvcenter{$\mapsto$} & 
    \myvcenter{%
    \begin{tikzpicture}
    [x={(-0.216506cm,-0.125000cm)}, y={(0.216506cm,-0.125000cm)}, z={(0.000000cm,0.250000cm)}]
    \draw[thick, densely dotted] (0,0,0) -- (1,0,0);
    \fill[fill=facecolor, draw=black, shift={(1,0,0)}]
    (0, 0, 0) -- (0, 0, 1) -- (1, 0, 1) -- (1, 0, 0) -- cycle;
    \node[circle,fill=black,draw=black,minimum size=1.2mm,inner sep=0pt] at (0,0,0) {};
    \end{tikzpicture}} \\
    \myvcenter{%
    \begin{tikzpicture}
    [x={(-0.216506cm,-0.125000cm)}, y={(0.216506cm,-0.125000cm)}, z={(0.000000cm,0.250000cm)}]
    \fill[fill=facecolor, draw=black, shift={(0,0,0)}]
    (0, 0, 0) -- (1, 0, 0) -- (1, 1, 0) -- (0, 1, 0) -- cycle;
    \node[circle,fill=black,draw=black,minimum size=1.2mm,inner sep=0pt] at (0,0,0) {};
    \end{tikzpicture}}%
     & \myvcenter{$\mapsto$} & 
    \myvcenter{%
    \begin{tikzpicture}
    [x={(-0.216506cm,-0.125000cm)}, y={(0.216506cm,-0.125000cm)}, z={(0.000000cm,0.250000cm)}]
    \draw[thick, densely dotted] (0,0,0) -- (1,0,0);
    \fill[fill=facecolor, draw=black, shift={(1,0,0)}]
    (0, 0, 0) -- (1, 0, 0) -- (1, 1, 0) -- (0, 1, 0) -- cycle;
    \node[circle,fill=black,draw=black,minimum size=1.2mm,inner sep=0pt] at (0,0,0) {};
    \end{tikzpicture}}
    \end{tabular}
\right.
\quad
\Sigma_2 : \left\{
    \begin{tabular}{rcl}%
    \myvcenter{%
    \begin{tikzpicture}
    [x={(-0.216506cm,-0.125000cm)}, y={(0.216506cm,-0.125000cm)}, z={(0.000000cm,0.250000cm)}]
    \fill[fill=facecolor, draw=black, shift={(0,0,0)}]
    (0, 0, 0) -- (0, 1, 0) -- (0, 1, 1) -- (0, 0, 1) -- cycle;
    \node[circle,fill=black,draw=black,minimum size=1.2mm,inner sep=0pt] at (0,0,0) {};
    \end{tikzpicture}}%
     & \myvcenter{$\mapsto$} & 
    \myvcenter{%
    \begin{tikzpicture}
    [x={(-0.216506cm,-0.125000cm)}, y={(0.216506cm,-0.125000cm)}, z={(0.000000cm,0.250000cm)}]
    \draw[thick, densely dotted] (0,0,0) -- (0,1,0);
    \fill[fill=facecolor, draw=black, shift={(0,1,0)}]
    (0, 0, 0) -- (0, 1, 0) -- (0, 1, 1) -- (0, 0, 1) -- cycle;
    \node[circle,fill=black,draw=black,minimum size=1.2mm,inner sep=0pt] at (0,0,0) {};
    \end{tikzpicture}} \\
    \myvcenter{%
    \begin{tikzpicture}
    [x={(-0.216506cm,-0.125000cm)}, y={(0.216506cm,-0.125000cm)}, z={(0.000000cm,0.250000cm)}]
    \fill[fill=facecolor, draw=black, shift={(0,0,0)}]
    (0, 0, 0) -- (0, 0, 1) -- (1, 0, 1) -- (1, 0, 0) -- cycle;
    \node[circle,fill=black,draw=black,minimum size=1.2mm,inner sep=0pt] at (0,0,0) {};
    \end{tikzpicture}}%
     & \myvcenter{$\mapsto$} & 
    \myvcenter{%
    \begin{tikzpicture}
    [x={(-0.216506cm,-0.125000cm)}, y={(0.216506cm,-0.125000cm)}, z={(0.000000cm,0.250000cm)}]
    \fill[fill=facecolor, draw=black, shift={(0,0,0)}]
    (0, 0, 0) -- (0, 1, 0) -- (0, 1, 1) -- (0, 0, 1) -- cycle;
    \fill[fill=facecolor, draw=black, shift={(0,0,0)}]
    (0, 0, 0) -- (0, 0, 1) -- (1, 0, 1) -- (1, 0, 0) -- cycle;
    \fill[fill=facecolor, draw=black, shift={(0,0,0)}]
    (0, 0, 0) -- (1, 0, 0) -- (1, 1, 0) -- (0, 1, 0) -- cycle;
    \node[circle,fill=black,draw=black,minimum size=1.2mm,inner sep=0pt] at (0,0,0) {};
    \end{tikzpicture}} \\
    \myvcenter{%
    \begin{tikzpicture}
    [x={(-0.216506cm,-0.125000cm)}, y={(0.216506cm,-0.125000cm)}, z={(0.000000cm,0.250000cm)}]
    \fill[fill=facecolor, draw=black, shift={(0,0,0)}]
    (0, 0, 0) -- (1, 0, 0) -- (1, 1, 0) -- (0, 1, 0) -- cycle;
    \node[circle,fill=black,draw=black,minimum size=1.2mm,inner sep=0pt] at (0,0,0) {};
    \end{tikzpicture}}%
     & \myvcenter{$\mapsto$} & 
    \myvcenter{%
    \begin{tikzpicture}
    [x={(-0.216506cm,-0.125000cm)}, y={(0.216506cm,-0.125000cm)}, z={(0.000000cm,0.250000cm)}]
    \draw[thick, densely dotted] (0,0,0) -- (0,1,0);
    \fill[fill=facecolor, draw=black, shift={(0,1,0)}]
    (0, 0, 0) -- (1, 0, 0) -- (1, 1, 0) -- (0, 1, 0) -- cycle;
    \node[circle,fill=black,draw=black,minimum size=1.2mm,inner sep=0pt] at (0,0,0) {};
    \end{tikzpicture}}
    \end{tabular}
\right.
\quad
\Sigma_3 \ : \ \left\{
    \begin{tabular}{rcl}%
    \myvcenter{%
    \begin{tikzpicture}
    [x={(-0.216506cm,-0.125000cm)}, y={(0.216506cm,-0.125000cm)}, z={(0.000000cm,0.250000cm)}]
    \fill[fill=facecolor, draw=black, shift={(0,0,0)}]
    (0, 0, 0) -- (0, 1, 0) -- (0, 1, 1) -- (0, 0, 1) -- cycle;
    \node[circle,fill=black,draw=black,minimum size=1.2mm,inner sep=0pt] at (0,0,0) {};
    \end{tikzpicture}}%
     & \myvcenter{$\mapsto$} & 
    \myvcenter{%
    \begin{tikzpicture}
    [x={(-0.216506cm,-0.125000cm)}, y={(0.216506cm,-0.125000cm)}, z={(0.000000cm,0.250000cm)}]
    \draw[thick, densely dotted] (0,0,0) -- (0,0,1);
    \fill[fill=facecolor, draw=black, shift={(0,0,1)}]
    (0, 0, 0) -- (0, 1, 0) -- (0, 1, 1) -- (0, 0, 1) -- cycle;
    \node[circle,fill=black,draw=black,minimum size=1.2mm,inner sep=0pt] at (0,0,0) {};
    \end{tikzpicture}} \\
    \myvcenter{%
    \begin{tikzpicture}
    [x={(-0.216506cm,-0.125000cm)}, y={(0.216506cm,-0.125000cm)}, z={(0.000000cm,0.250000cm)}]
    \fill[fill=facecolor, draw=black, shift={(0,0,0)}]
    (0, 0, 0) -- (0, 0, 1) -- (1, 0, 1) -- (1, 0, 0) -- cycle;
    \node[circle,fill=black,draw=black,minimum size=1.2mm,inner sep=0pt] at (0,0,0) {};
    \end{tikzpicture}}%
     & \myvcenter{$\mapsto$} & 
    \myvcenter{%
    \begin{tikzpicture}
    [x={(-0.216506cm,-0.125000cm)}, y={(0.216506cm,-0.125000cm)}, z={(0.000000cm,0.250000cm)}]
    \draw[thick, densely dotted] (0,0,0) -- (0,0,1);
    \fill[fill=facecolor, draw=black, shift={(0,0,1)}]
    (0, 0, 0) -- (0, 0, 1) -- (1, 0, 1) -- (1, 0, 0) -- cycle;
    \node[circle,fill=black,draw=black,minimum size=1.2mm,inner sep=0pt] at (0,0,0) {};
    \end{tikzpicture}} \\
    \myvcenter{%
    \begin{tikzpicture}
    [x={(-0.216506cm,-0.125000cm)}, y={(0.216506cm,-0.125000cm)}, z={(0.000000cm,0.250000cm)}]
    \fill[fill=facecolor, draw=black, shift={(0,0,0)}]
    (0, 0, 0) -- (1, 0, 0) -- (1, 1, 0) -- (0, 1, 0) -- cycle;
    \node[circle,fill=black,draw=black,minimum size=1.2mm,inner sep=0pt] at (0,0,0) {};
    \end{tikzpicture}}%
     & \myvcenter{$\mapsto$} & 
    \myvcenter{%
    \begin{tikzpicture}
    [x={(-0.216506cm,-0.125000cm)}, y={(0.216506cm,-0.125000cm)}, z={(0.000000cm,0.250000cm)}]
    \fill[fill=facecolor, draw=black, shift={(0,0,0)}]
    (0, 0, 0) -- (0, 1, 0) -- (0, 1, 1) -- (0, 0, 1) -- cycle;
    \fill[fill=facecolor, draw=black, shift={(0,0,0)}]
    (0, 0, 0) -- (0, 0, 1) -- (1, 0, 1) -- (1, 0, 0) -- cycle;
    \fill[fill=facecolor, draw=black, shift={(0,0,0)}]
    (0, 0, 0) -- (1, 0, 0) -- (1, 1, 0) -- (0, 1, 0) -- cycle;
    \node[circle,fill=black,draw=black,minimum size=1.2mm,inner sep=0pt] at (0,0,0) {};
    \end{tikzpicture}}
    \end{tabular}
\right..
\]

\begin{exam}
Below are some examples of images of some patterns by
the dual substitutions $\Sigma_1$, $\Sigma_2$, $\Sigma_3$.
Let $\mcU = [\mathbf 0,1]^* \cup [\mathbf 0,2]^* \cup [\mathbf 0,3]^*$.
\[
\renewcommand{\arraystretch}{1.5}
\begin{array}{rclcl}
\Sigma_1(\mcU) & = &
\mcU \cup [(1, 0, 0), 2]^* \cup [(1, 0, 0), 3]^* & = &
\myvcenter{%
\begin{tikzpicture}
[x={(-0.216506cm,-0.125000cm)}, y={(0.216506cm,-0.125000cm)}, z={(0.000000cm,0.250000cm)}]
\definecolor{facecolor}{rgb}{0.8,0.8,0.8}
\fill[fill=facecolor, draw=black, shift={(0,0,0)}]
(0, 0, 0) -- (0, 0, 1) -- (1, 0, 1) -- (1, 0, 0) -- cycle;
\fill[fill=facecolor, draw=black, shift={(0,0,0)}]
(0, 0, 0) -- (0, 1, 0) -- (0, 1, 1) -- (0, 0, 1) -- cycle;
\fill[fill=facecolor, draw=black, shift={(1,0,0)}]
(0, 0, 0) -- (1, 0, 0) -- (1, 1, 0) -- (0, 1, 0) -- cycle;
\fill[fill=facecolor, draw=black, shift={(1,0,0)}]
(0, 0, 0) -- (0, 0, 1) -- (1, 0, 1) -- (1, 0, 0) -- cycle;
\fill[fill=facecolor, draw=black, shift={(0,0,0)}]
(0, 0, 0) -- (1, 0, 0) -- (1, 1, 0) -- (0, 1, 0) -- cycle;
\node[circle,fill=black,draw=black,minimum size=1.2mm,inner sep=0pt] at (0,0,0) {};
\end{tikzpicture}} \\
%%%
\Sigma_2(\mcU) & = &
\mcU \cup [(0, 1, 0), 1]^* \cup [(0, 1, 0), 3]^* & = &
\myvcenter{%
\begin{tikzpicture}
[x={(-0.216506cm,-0.125000cm)}, y={(0.216506cm,-0.125000cm)}, z={(0.000000cm,0.250000cm)}]
\definecolor{facecolor}{rgb}{0.8,0.8,0.8}
\fill[fill=facecolor, draw=black, shift={(0,0,0)}]
(0, 0, 0) -- (0, 0, 1) -- (1, 0, 1) -- (1, 0, 0) -- cycle;
\fill[fill=facecolor, draw=black, shift={(0,0,0)}]
(0, 0, 0) -- (0, 1, 0) -- (0, 1, 1) -- (0, 0, 1) -- cycle;
\fill[fill=facecolor, draw=black, shift={(0,1,0)}]
(0, 0, 0) -- (1, 0, 0) -- (1, 1, 0) -- (0, 1, 0) -- cycle;
\fill[fill=facecolor, draw=black, shift={(0,0,0)}]
(0, 0, 0) -- (1, 0, 0) -- (1, 1, 0) -- (0, 1, 0) -- cycle;
\fill[fill=facecolor, draw=black, shift={(0,1,0)}]
(0, 0, 0) -- (0, 1, 0) -- (0, 1, 1) -- (0, 0, 1) -- cycle;
\node[circle,fill=black,draw=black,minimum size=1.2mm,inner sep=0pt] at (0,0,0) {};
\end{tikzpicture}} \\
%%%
\Sigma_3(\mcU) & = &
\mcU \cup [(0, 0, 1), 1]^* \cup [(0, 0, 1), 2]^* & = &
\myvcenter{%
\begin{tikzpicture}
[x={(-0.216506cm,-0.125000cm)}, y={(0.216506cm,-0.125000cm)}, z={(0.000000cm,0.250000cm)}]
\definecolor{facecolor}{rgb}{0.8,0.8,0.8}
\fill[fill=facecolor, draw=black, shift={(0,0,1)}]
(0, 0, 0) -- (0, 0, 1) -- (1, 0, 1) -- (1, 0, 0) -- cycle;
\fill[fill=facecolor, draw=black, shift={(0,0,0)}]
(0, 0, 0) -- (0, 0, 1) -- (1, 0, 1) -- (1, 0, 0) -- cycle;
\fill[fill=facecolor, draw=black, shift={(0,0,0)}]
(0, 0, 0) -- (0, 1, 0) -- (0, 1, 1) -- (0, 0, 1) -- cycle;
\fill[fill=facecolor, draw=black, shift={(0,0,1)}]
(0, 0, 0) -- (0, 1, 0) -- (0, 1, 1) -- (0, 0, 1) -- cycle;
\fill[fill=facecolor, draw=black, shift={(0,0,0)}]
(0, 0, 0) -- (1, 0, 0) -- (1, 1, 0) -- (0, 1, 0) -- cycle;
\node[circle,fill=black,draw=black,minimum size=1.2mm,inner sep=0pt] at (0,0,0) {};
\end{tikzpicture}} \\
%%%
\Sigma_2(\Sigma_1(\mcU)) & = &
\Sigma_2(\mcU) \cup [(1, 0, 0), 3]^* \cup (\mcU + (1,-1,0)) & = &
\myvcenter{%
\begin{tikzpicture}
[x={(-0.216506cm,-0.125000cm)}, y={(0.216506cm,-0.125000cm)}, z={(0.000000cm,0.250000cm)}]
\definecolor{facecolor}{rgb}{0.8,0.8,0.8}
\fill[fill=facecolor, draw=black, shift={(1,0,0)}]
(0, 0, 0) -- (1, 0, 0) -- (1, 1, 0) -- (0, 1, 0) -- cycle;
\fill[fill=facecolor, draw=black, shift={(1,-1,0)}]
(0, 0, 0) -- (0, 0, 1) -- (1, 0, 1) -- (1, 0, 0) -- cycle;
\fill[fill=facecolor, draw=black, shift={(0,1,0)}]
(0, 0, 0) -- (0, 1, 0) -- (0, 1, 1) -- (0, 0, 1) -- cycle;
\fill[fill=facecolor, draw=black, shift={(0,1,0)}]
(0, 0, 0) -- (1, 0, 0) -- (1, 1, 0) -- (0, 1, 0) -- cycle;
\fill[fill=facecolor, draw=black, shift={(0,0,0)}]
(0, 0, 0) -- (0, 0, 1) -- (1, 0, 1) -- (1, 0, 0) -- cycle;
\fill[fill=facecolor, draw=black, shift={(1,-1,0)}]
(0, 0, 0) -- (0, 1, 0) -- (0, 1, 1) -- (0, 0, 1) -- cycle;
\fill[fill=facecolor, draw=black, shift={(0,0,0)}]
(0, 0, 0) -- (1, 0, 0) -- (1, 1, 0) -- (0, 1, 0) -- cycle;
\fill[fill=facecolor, draw=black, shift={(1,-1,0)}]
(0, 0, 0) -- (1, 0, 0) -- (1, 1, 0) -- (0, 1, 0) -- cycle;
\fill[fill=facecolor, draw=black, shift={(0,0,0)}]
(0, 0, 0) -- (0, 1, 0) -- (0, 1, 1) -- (0, 0, 1) -- cycle;
\node[circle,fill=black,draw=black,minimum size=1.2mm,inner sep=0pt] at (0,0,0) {};
\end{tikzpicture}}
\end{array}
\]
\end{exam}

In the following we will also need to ``desubstitute'', \emph{i.e.}, to compute preimages of faces.
The following lemma gives the preimages of unit faces for the Arnoux-Rauzy dual substitutions.
The proof, which is an easy enumeration of cases, can be found in Appendix \ref{app:preimages}.
By abuse of notation, $\Sigma_k^{-1}([\bfx, i]^*)$ stands for $\Sigma_k^{-1}(\{[\bfx, i]^*\})$,
\emph{i.e.}, for the union of faces $[\bfy, j]^*$ such that $\Sigma_k([\bfy, j]^*)$ contains the face $[\bfx, i]^*$.
This notation will also be extended similarly to finite unions of faces in all that follows.
\begin{lemm}[Preimages of faces]
\label{lemm:preimages}
Let $\bfx =\svect{x}{y}{z} \in \bbZ^3$.
We have
\[
\Sigma_1^{-1}([\bfx, i]^*) \ = \
\left\{
\begin{array}{ll}
  \big[ \svect{x+y+z}{y}{z}, 1 \big]^*
    & \textup{if } i = 1 \\
  \big[ \svect{x+y+z}{y}{z}, 1 \big]^* \cup \big[ \svect{x+y+z-1}{y}{z}, i \big]^*
    & \textup{if } i = 2,3
\end{array}
\right.,
\]
\[
\Sigma_2^{-1}([\bfx, i]^*) \ = \
\left\{
\begin{array}{ll}
  \big[ \svect{x}{x+y+z}{z}, 2 \big]^*
    & \textup{if } i = 2 \\
  \big[ \svect{x}{x+y+z}{z}, 2 \big]^* \cup \big[ \svect{x}{x+y+z-1}{z}, i \big]^*
    & \textup{if } i = 1,3
\end{array}
\right.,
\]
\[
\Sigma_3^{-1}([\bfx, i]^*) \ = \
\left\{
\begin{array}{ll}
  \big[ \svect{x}{y}{x+y+z}, 3 \big]^*
    & \textup{if } i = 3 \\
  \big[ \svect{x}{y}{x+y+z}, 3 \big]^* \cup \big[ \svect{x}{y}{x+y+z-1}, i \big]^*
    & \textup{if } i = 2,3
\end{array}
\right..
\]
\end{lemm}

\subsection{Rauzy fractals and spectrum}
\label{subsec:RF}
\paragraph{Geometrical representations and spectra}
A fundamental question concerning substitutive dynamical systems $(X_\sigma, S)$
deals with the possibility of giving them a geometric representation.
By a geometric representation, we mean a dynamical system of geometric nature
that is topologically conjugate or measure-theoretically isomorphic to $(X_\sigma, S)$.
In particular, we are interested in conditions under which
it is possible to give a geometric representation of $(X_\sigma, S)$
as a translation on the torus (or more generally on a locally compact abelian group).

This question can be reformulated in spectral terms.
Indeed, one associates with the (uniquely ergodic) symbolic dynamical system $(X_\sigma, S)$
(considered as a measure-theoretic dynamical system)
the operator $U$ acting on the Hilbert space $\mcL^{2}(X,\mu)$
(where $\mu$ stands for the invariant measure) defined as follows:
\[
U \ : \ \left\{
\begin{array}{rcl}
    \mcL^2(X, \mu) & \rightarrow & \mcL^2(X, \mu) \\
     f & \mapsto & f \circ T
\end{array}
\right..
\]
The \emph{eigenvalues} and \emph{eigenfunctions} of $(X, T,\mu)$
are defined as being those of the unitary operator $U$.
The symbolic dynamical system $(X_\sigma, S)$ is said to have \emph{pure discrete spectrum}
if $\mcL^2(X, \mu)$ admits an Hilbert basis of eigenfunctions.
According to Von Neumann's theorem,
$(X_\sigma, S)$ has pure discrete spectrum if and only if
it is conjugate to a translation on a compact group.
For more details, see \cite{Wal82, Que10}.

The question of the existence of a geometric representation can thus be rephrased as follows:
which are the substitutions whose associated dynamical system has discrete spectrum?
Recall that measure-theoretic discrete spectrum and topological
discrete spectrum are proved to be equivalent for symbolic dynamical systems
associated with a primitive substitution \cite{Hos86}.

\paragraph{Rauzy fractals}
It is widely believed that every Pisot irreducible substitution has pure discrete spectrum;
this is known as the \emph{Pisot conjecture}
(see \cite[Chap. 7]{PF02} and \cite{BK06} for more details).
A strategy for providing geometric representations (and tackling the Pisot conjecture)
has been initiated by Rauzy in \cite{Rau82}
for the particular Tribonacci substitution
$\sigma : 1 \mapsto 12, \ 2 \mapsto 13, \ 3 \mapsto 1$,
where it is proved that $(X_\sigma, S)$ is measure-theoretically isomorphic
to the two dimensional toral translation given by $x \mapsto x + (1/\beta,1/\beta^2)$
($\beta$ is the Pisot eigenvalue of $\Ms$).

\emph{Rauzy fractals} were first introduced in \cite{Rau82} in the case of the Tribonacci substitution,
and then in \cite{Thu89} in the case of the $\beta$-numeration associated with the smallest Pisot number.
They are good candidates to be fundamental domains for the toral translation
that is conjectured to be measure-theoretically isomorphic to $(X_\sigma, S)$.
Such fractals can more generally be associated with Pisot substitutions (see \cite{BS05, ST10, PF02}),
as well as with Pisot $\beta$-shifts under the name of central tiles (see \cite{Aki98, Aki02}).

An effective combinatorial construction of Rauzy fractals,
different from the original one, has been proposed in \cite{AI01}.
Let $\sigma $ be a unimodular Pisot irreducible substitution.
Remark~\ref{rema:imgU} and Proposition~\ref{prop:imgplane} enable us to iterate $\EOSS$
on $\mcU = [\mathbf 0,1]^* \cup [\mathbf 0,2]^* \cup [\mathbf 0,3]^*$
in order to obtain an infinite increasing sequence of patterns (with respect to inclusion).
It is then possible to project and renormalize the patterns $\EOSS^n(\mcU)$
by applying $\Ms \circ \pic$ (see Definition \ref{defi:contrac}),
in order to obtain a sequence of subsets of the contracting plane $\Pc$
that converges to a compact subset of $\Pc$.
More precisely, for $n \geq 0$,
let $(\mcD_n)_{n \geq 0}$ be the sequence of \emph{approximations}
of the Rauzy fractal associated with $\sigma$, defined for all $n$ by
\[
\mcD_n \ = \ \Ms^n \circ \pic \circ \EOSS^n(\mcU).
\]
Arnoux and Ito proved in~\cite{AI01} that $(\mcD_n)_{n \geq 0}$ is a convergent sequence
in the metric space of compact subsets of $\Pc$ together with the Hausdorff distance.
This yields the following definition for the Rauzy fractal
associated with a substitution.

\begin{defi}[Rauzy fractal]
\label{def:RF}
Let $\sigma $ be a unimodular Pisot irreducible substitution.
The \emph{Rauzy fractal} associated with $\sigma$ is the Hausdorff limit
of the sequence of approximations $(\mcD_n)_{n \geq 0}$.
\end{defi}

This definition allows us to give an equivalent formulation of the pure discrete spectrum property
in terms of the approximations $\mcD_n$ of the Rauzy fractal.
The ``if'' direction comes from \cite{AI01, IR05}
and the ``only if'' direction has been proved in \cite{BK06}
(see also \cite{BST11}).

\begin{theo}
\label{theo:puredis}
Let $\sigma$ be a unimodular Pisot irreducible substitution.
The substitutive symbolic dynamical system $(X_\sigma,S)$ has pure discrete spectrum if and only if
the approximation set $\mcD_n$ contains arbitrarily large balls, as $n$ goes to infinity.
\end{theo}

Our aim is now to prove that the approximations $ \mcD_n$ contain arbitrarily large balls when $n$ tends to infinity,
when $\sigma$ is a product of Arnoux-Rauzy substitutions.
This will be the object of Sections \ref{sec:cov} and \ref{subsec:annulus}.

\section{Covering properties and annuli}
\label{sec:cov}
We now introduce the combinatorial tools that will be used in Section \ref{subsec:annulus},
in order to prove that some increasing sequences of patterns contain arbitrarily large balls.

A finite union of unit faces is called a \emph{pattern}.
If $\mcL$ is a set of patterns, we denote by
$\omcL = \{v + P : P \in \mcL, v \in \bbZ^3 \}$
the set of all the translations of patterns of $\mcL$.
In the following, we will not make the distinction between a pattern and the set of its faces.

\subsection{Coverings}
\label{subsec:cov}
The covering property below (introduced in \cite{IO93})
can be seen as an arc-connectedness type property with respect to a set of patterns $\mcL$:
every couple of faces of an $\mcL$-covered pattern
must be connected by a path of overlapping patterns of $\mcL$.
If $\mcL$ is well chosen, $\mcL$-covering guarantees good topological properties for the sets $\mcD_n$,
such as connectedness (see, \emph{e.g.}, \cite{BJS11} for more details).

\begin{defi}[$\mcL$-cover]
\label{defi:cover}
Let $\mcL$ be a set of patterns.
A pattern $P$ is \emph{$\mcL$-covered} if for all faces $e, f \in P$,
there exists $(Q_1, \ldots, Q_n) \in \omcL^n$ such that:
\begin{enumerate}
  \item $e \in Q_1$ and $f \in Q_n$;
  \item $Q_k \cap Q_{k+1}$ contains at least one face,
    for all $k \in \{1, \ldots, n-1\}$; 
  \item $Q_k \subseteq P$ for all $k \in \{1, \ldots, n\}$.
\end{enumerate}
\end{defi}

The following proposition gives a convenient combinatorial condition to
prove that $\mcL$-covering of the images of a pattern by a dual substitution.

\begin{prop}[\cite{IO93}]
\label{prop:coverprop}
Let $P$ be an $\mcL$-covered pattern,
$\Sigma$ be a dual substitution
and $\mcL$ be a set of patterns such that
$\Sigma(Q)$ is $\mcL$-covered for all $Q \in \mcL$.
Then $\Sigma(P)$ is $\mcL$-covered.
\end{prop}

\subsection{Strong coverings}
The notion of $\mcL$-covering defined in Section \ref{subsec:cov}
is already enough to provide us with some of the topological properties
which will be needed in Section~\ref{subsec:annulus}.
However, to ensure that these properties are preserved under the image of a dual substitution,
we need a stronger property, namely \emph{strong} $\mcL$-covering,
that we define in this section and that is the main new ingredient of this article
(this notion plays a key role in particular in the proof of Lemma~\ref{lemm:annulus_induction}).
Example \ref{exam:broken_annulus} shows a concrete situation where strong $\mcL$-covers are needed (see also Remark~\ref{rema:IO}).

Before defining strong $\mcL$-covers,
we introduce the set of patterns $\Ledge$
that consists of all the twelve edge-connected two-face patterns defined below
(a pattern is \emph{edge-connected} if it is $\Ledge$-covered):
\[
\renewcommand{\arraystretch}{1.3}
\begin{array}{cclllccclll}
\input{fig/E11a.tex} & = & [\mathbf 0, 1]^* & \cup & [(0,1,0), 1]^* & \ &
\input{fig/E12a.tex} & = & [\mathbf 0, 1]^* & \cup & [\mathbf 0, 2]^* \\
\input{fig/E11b.tex} & = & [\mathbf 0, 1]^* & \cup & [(0,0,1), 1]^* & \ &
\input{fig/E12b.tex} & = & [\mathbf 0, 1]^* & \cup & [(-1,1,0), 2]^* \\
\input{fig/E22a.tex} & = & [\mathbf 0, 2]^* & \cup & [(1,0,0), 2]^* & \ &
\input{fig/E13a.tex} & = & [\mathbf 0, 1]^* & \cup & [\mathbf 0, 3]^* \\
\input{fig/E22b.tex} & = & [\mathbf 0, 2]^* & \cup & [(0,0,1), 2]^* & \ &
\input{fig/E13b.tex} & = & [\mathbf 0, 1]^* & \cup & [(-1,0,1), 3]^* \\
\input{fig/E33a.tex} & = & [\mathbf 0, 3]^* & \cup & [(1,0,0), 3]^* & \ &
\input{fig/E23a.tex} & = & [\mathbf 0, 2]^* & \cup & [\mathbf 0, 3]^* \\
\input{fig/E33b.tex} & = & [\mathbf 0, 3]^* & \cup & [(0,1,0), 3]^* & \ &
\input{fig/E23b.tex} & = & [\mathbf 0, 2]^* & \cup & [(0,-1,1), 3]^*.
\end{array}
\]

\begin{defi}[Strong $\mcL$-cover]
Let $\mcL$ be a set of edge-connected patterns.
A pattern $P$ is \emph{strongly $\mcL$-covered} if
\begin{enumerate}
\item $P$ is $\mcL$-covered;
\item for every pattern $X \in \oLedge$ such that $X \subseteq P$,
    there exists a pattern $Y \in \omcL$ such that $X \subseteq Y \subseteq P$.
\end{enumerate}
\end{defi}

The intuitive idea behind the notion of strong $\mcL$-covering is that every occurrence
of a pattern of $\oLedge$ in $P$ is required to be ``completed within $P$'' by a pattern of $\mcL$.

\begin{exam}
Let
$P =
\myvcenter{\begin{tikzpicture}
[x={(-0.216506cm,-0.125000cm)}, y={(0.216506cm,-0.125000cm)}, z={(0.000000cm,0.250000cm)}]
\definecolor{facecolor}{rgb}{0.800,0.800,0.800}
\fill[fill=facecolor, draw=black, shift={(0,0,0)}]
(0, 0, 0) -- (0, 1, 0) -- (0, 1, 1) -- (0, 0, 1) -- cycle;
\fill[fill=facecolor, draw=black, shift={(0,0,0)}]
(0, 0, 0) -- (0, 0, 1) -- (1, 0, 1) -- (1, 0, 0) -- cycle;
\fill[fill=facecolor, draw=black, shift={(0,0,0)}]
(0, 0, 0) -- (1, 0, 0) -- (1, 1, 0) -- (0, 1, 0) -- cycle;
\end{tikzpicture}}\,$ and
$Q =
\myvcenter{\begin{tikzpicture}
[x={(-0.216506cm,-0.125000cm)}, y={(0.216506cm,-0.125000cm)}, z={(0.000000cm,0.250000cm)}]
\definecolor{facecolor}{rgb}{0.800,0.800,0.800}
\fill[fill=facecolor, draw=black, shift={(0,0,0)}]
(0, 0, 0) -- (0, 1, 0) -- (0, 1, 1) -- (0, 0, 1) -- cycle;
\fill[fill=facecolor, draw=black, shift={(0,0,0)}]
(0, 0, 0) -- (0, 0, 1) -- (1, 0, 1) -- (1, 0, 0) -- cycle;
\fill[fill=facecolor, draw=black, shift={(0,0,0)}]
(0, 0, 0) -- (1, 0, 0) -- (1, 1, 0) -- (0, 1, 0) -- cycle;
\fill[fill=facecolor, draw=black, shift={(1,0,0)}]
(0, 0, 0) -- (1, 0, 0) -- (1, 1, 0) -- (0, 1, 0) -- cycle;
\end{tikzpicture}}\,$,
and let
$\mcL = \big\{\,
\input{fig/E12a_nodot.tex},
\input{fig/E13a_nodot.tex},
\input{fig/E33aAR_nodot.tex}
\big\}$.
Both $P$ and $Q$ are $\mcL$-covered,
but only $Q$ is strongly $\mcL$-covered
because the two-face pattern $X = \input{fig/E23a_nodot.tex}$
can be ``completed'' by \input{fig/E33aAR_nodot.tex} within $Q$ but not within $P$.
Another example of an $\mcL$-covered pattern which is not strongly $\mcL$-covered
appears in Remark \ref{rema:contreexemplestrongcovering}.
\end{exam}

In the above definition of strong $\mcL$-covers,
we require $X$ to be in $\oLedge$
because it is precisely what we need for the results of this article,
due to technical reasons based on the particular form of the Arnoux-Rauzy substitutions,
but the same kind of definition could be made with another set of patterns than $\oLedge$.

Unfortunately, an analog of Proposition \ref{prop:coverprop}
does not hold in general for \emph{strong} $\mcL$-covers,
but we will prove an \emph{ad hoc} version of it in the particular case of
Arnoux-Rauzy substitutions (see Proposition \ref{prop:strongcovAR} below).

\subsection{Annuli}
The aim of this section is to define a notion of topological annulus
in the setting of finite unions of unit faces.
The main argument of this paper is based on this object:
we will prove in Section \ref{subsec:annulus} that,
under some technical hypotheses,
the image of an annulus by a dual substitution remains an annulus,
thanks to an \emph{ad hoc} generalization of Proposition~\ref{prop:coverprop}
(namely, Proposition \ref{prop:strongcovAR}).
This will enable us to deduce in Section \ref{subsec:dynconseq} that
the approximations $\mcD_n$ of Rauzy fractals contain arbitrarily large balls when $n$ tends to infinity.

\begin{defi}[Annulus]
\label{defi:annulus}
Let $\mcL$ be a set of edge-connected patterns and $\Gamma$ be a discrete plane.
An \emph{$\mcL$-annulus} of a pattern $P \subseteq \Gamma$
is a pattern $A \subseteq \Gamma$ such that:
\begin{enumerate}
  \item $P$, $A \cup P$ and $\Gamma \setminus (A \cup P)$ are $\mcL$-covered;
    \label{defi:anneauprop1}
  \item $A$ is strongly $\mcL$-covered;
    \label{defi:anneauprop2}
  \item $A$ and $P$ have no face in common;
    \label{defi:anneauprop3}
  \item $P \cap \overline{\Gamma \setminus (P \cup A)} = \varnothing$.
    \label{defi:anneauprop4}
\end{enumerate}
\end{defi}
Topological features of annuli are guarantied by 
Conditions (\ref{defi:anneauprop3}) and (\ref{defi:anneauprop4}).
The intuitive meaning of Condition (\ref{defi:anneauprop4}) is that $A$ completely ``surrounds'' $P$,
in such that a way that $P$ does not touch the ``outside'' region $\overline{\Gamma \setminus (P \cup A)}$.
Conditions (\ref{defi:anneauprop1}) and (\ref{defi:anneauprop2})
provide a convenient combinatorial way to ensure that
the topological features of Conditions (\ref{defi:anneauprop3}) and (\ref{defi:anneauprop4})
are preserved under the image by Arnoux-Rauzy substitutions,
as will be proved in Section \ref{subsec:annulus}.

\begin{rema}\label{rema:IO}
Contrarily to what is claimed in \cite{IO94},
the image by a dual substitution of a pattern that satisfies Conditions (\ref{defi:anneauprop1}),(\ref{defi:anneauprop3}) and (\ref{defi:anneauprop4})
has no reason to satisfy again these conditions,
as illustrated by Example \ref{exam:broken_annulus} below.
Nevertheless, the notion of strong covering with respect to a well-chosen set of patterns
will allow us in the case of the Arnoux-Rauzy substitutions to prove that
the image of an annulus remains an annulus.
\end{rema} 

\begin{exam}
Let $\LAR$ be the set of patterns defined in Section \ref{subsec:coverAR}, and let
\[
P_1 \cup A_1 \ = \
\myvcenter{%
\begin{tikzpicture}
[x={(-0.216506cm,-0.125000cm)}, y={(0.216506cm,-0.125000cm)}, z={(0.000000cm,0.250000cm)}]
\definecolor{facecolor}{rgb}{0.800,0.800,0.800}
\fill[fill=facecolor, draw=black, shift={(-1,1,1)}]
(0, 0, 0) -- (1, 0, 0) -- (1, 1, 0) -- (0, 1, 0) -- cycle;
\fill[fill=facecolor, draw=black, shift={(1,0,0)}]
(0, 0, 0) -- (1, 0, 0) -- (1, 1, 0) -- (0, 1, 0) -- cycle;
\fill[fill=facecolor, draw=black, shift={(-1,0,1)}]
(0, 0, 0) -- (1, 0, 0) -- (1, 1, 0) -- (0, 1, 0) -- cycle;
\fill[fill=facecolor, draw=black, shift={(0,1,0)}]
(0, 0, 0) -- (0, 1, 0) -- (0, 1, 1) -- (0, 0, 1) -- cycle;
\fill[fill=facecolor, draw=black, shift={(-1,0,1)}]
(0, 0, 0) -- (0, 0, 1) -- (1, 0, 1) -- (1, 0, 0) -- cycle;
\fill[fill=facecolor, draw=black, shift={(0,1,0)}]
(0, 0, 0) -- (1, 0, 0) -- (1, 1, 0) -- (0, 1, 0) -- cycle;
\definecolor{facecolor}{rgb}{0.350,0.350,0.350}
\fill[fill=facecolor, draw=black, shift={(0,0,0)}]
(0, 0, 0) -- (0, 0, 1) -- (1, 0, 1) -- (1, 0, 0) -- cycle;
\definecolor{facecolor}{rgb}{0.800,0.800,0.800}
\fill[fill=facecolor, draw=black, shift={(0,-1,1)}]
(0, 0, 0) -- (1, 0, 0) -- (1, 1, 0) -- (0, 1, 0) -- cycle;
\fill[fill=facecolor, draw=black, shift={(1,-1,0)}]
(0, 0, 0) -- (0, 1, 0) -- (0, 1, 1) -- (0, 0, 1) -- cycle;
\definecolor{facecolor}{rgb}{0.350,0.350,0.350}
\fill[fill=facecolor, draw=black, shift={(0,0,0)}]
(0, 0, 0) -- (1, 0, 0) -- (1, 1, 0) -- (0, 1, 0) -- cycle;
\definecolor{facecolor}{rgb}{0.800,0.800,0.800}
\fill[fill=facecolor, draw=black, shift={(1,-1,0)}]
(0, 0, 0) -- (1, 0, 0) -- (1, 1, 0) -- (0, 1, 0) -- cycle;
\fill[fill=facecolor, draw=black, shift={(0,-1,1)}]
(0, 0, 0) -- (0, 1, 0) -- (0, 1, 1) -- (0, 0, 1) -- cycle;
\definecolor{facecolor}{rgb}{0.350,0.350,0.350}
\fill[fill=facecolor, draw=black, shift={(0,0,0)}]
(0, 0, 0) -- (0, 1, 0) -- (0, 1, 1) -- (0, 0, 1) -- cycle;
\end{tikzpicture}%
}
\qquad \text{and} \qquad
P_2 \cup A_2 \ = \
\myvcenter{%
\begin{tikzpicture}
[x={(-0.216506cm,-0.125000cm)}, y={(0.216506cm,-0.125000cm)}, z={(0.000000cm,0.250000cm)}]
\definecolor{facecolor}{rgb}{0.800,0.800,0.800}
\fill[fill=facecolor, draw=black, shift={(0,0,1)}]
(0, 0, 0) -- (0, 0, 1) -- (1, 0, 1) -- (1, 0, 0) -- cycle;
\definecolor{facecolor}{rgb}{0.350,0.350,0.350}
\fill[fill=facecolor, draw=black, shift={(0,0,0)}]
(0, 0, 0) -- (0, 1, 0) -- (0, 1, 1) -- (0, 0, 1) -- cycle;
\definecolor{facecolor}{rgb}{0.800,0.800,0.800}
\fill[fill=facecolor, draw=black, shift={(0,0,1)}]
(0, 0, 0) -- (0, 1, 0) -- (0, 1, 1) -- (0, 0, 1) -- cycle;
\fill[fill=facecolor, draw=black, shift={(0,1,0)}]
(0, 0, 0) -- (0, 1, 0) -- (0, 1, 1) -- (0, 0, 1) -- cycle;
\fill[fill=facecolor, draw=black, shift={(-1,1,1)}]
(0, 0, 0) -- (1, 0, 0) -- (1, 1, 0) -- (0, 1, 0) -- cycle;
\fill[fill=facecolor, draw=black, shift={(1,-1,1)}]
(0, 0, 0) -- (0, 1, 0) -- (0, 1, 1) -- (0, 0, 1) -- cycle;
\fill[fill=facecolor, draw=black, shift={(1,0,-1)}]
(0, 0, 0) -- (0, 1, 0) -- (0, 1, 1) -- (0, 0, 1) -- cycle;
\definecolor{facecolor}{rgb}{0.350,0.350,0.350}
\fill[fill=facecolor, draw=black, shift={(0,0,0)}]
(0, 0, 0) -- (0, 0, 1) -- (1, 0, 1) -- (1, 0, 0) -- cycle;
\definecolor{facecolor}{rgb}{0.800,0.800,0.800}
\fill[fill=facecolor, draw=black, shift={(0,1,-1)}]
(0, 0, 0) -- (0, 0, 1) -- (1, 0, 1) -- (1, 0, 0) -- cycle;
\fill[fill=facecolor, draw=black, shift={(1,-1,0)}]
(0, 0, 0) -- (0, 1, 0) -- (0, 1, 1) -- (0, 0, 1) -- cycle;
\fill[fill=facecolor, draw=black, shift={(0,1,-1)}]
(0, 0, 0) -- (0, 1, 0) -- (0, 1, 1) -- (0, 0, 1) -- cycle;
\definecolor{facecolor}{rgb}{0.350,0.350,0.350}
\fill[fill=facecolor, draw=black, shift={(0,0,0)}]
(0, 0, 0) -- (1, 0, 0) -- (1, 1, 0) -- (0, 1, 0) -- cycle;
\definecolor{facecolor}{rgb}{0.800,0.800,0.800}
\fill[fill=facecolor, draw=black, shift={(1,-1,0)}]
(0, 0, 0) -- (1, 0, 0) -- (1, 1, 0) -- (0, 1, 0) -- cycle;
\fill[fill=facecolor, draw=black, shift={(1,-1,0)}]
(0, 0, 0) -- (0, 0, 1) -- (1, 0, 1) -- (1, 0, 0) -- cycle;
\fill[fill=facecolor, draw=black, shift={(1,0,-1)}]
(0, 0, 0) -- (0, 0, 1) -- (1, 0, 1) -- (1, 0, 0) -- cycle;
\fill[fill=facecolor, draw=black, shift={(-1,1,1)}]
(0, 0, 0) -- (0, 0, 1) -- (1, 0, 1) -- (1, 0, 0) -- cycle;
\end{tikzpicture}%
}
\]
where $P_1 = P_2 = \mcU$ is shown in dark gray and $A_1$, $A_2$ are in light gray.
(The underlying discrete plane $\Gamma$ is not shown.)
The pattern $A_1$ is not an $\LAR$-annulus of $P_1$,
because Condition \ref{defi:anneauprop4} fails
(the lower vertex of $\mcU$ intersects $\overline{\Gamma \setminus (P_1 \cup A_1)}$).
However, $A_2$ is an $\LAR$-annulus of $P_2$:
Conditions (\ref{defi:anneauprop1}) and (\ref{defi:anneauprop2}) are easy to check
and Conditions (\ref{defi:anneauprop3}) and (\ref{defi:anneauprop4}) clearly hold.
\end{exam}

\begin{exam}
\label{exam:broken_annulus}
Let $P \cup A \cup f$ be the pattern shown below.
\begin{equation}
\label{eq:compl}
P \cup A \cup f \quad = \quad
\input{fig/broken_annulus1.tex}
\qquad \stackrel{\Sigma_2\Sigma_2}{\longmapsto} \qquad
\input{fig/broken_annulus2.tex}
\end{equation}
($P$ is shown in dark gray, $A$ in light gray and $f$ in white.)
The pattern $A$ is not an $\LAR$-annulus of $P$,
because it is not strongly $\LAR$-covered.
The strong $\LAR$-covering fails because of the two-face pattern of $A$ shown below in white:
\[
\input{fig/broken_annulus1_bis.tex}.
\]
We can observe that $\Sigma_2 \Sigma_2(A)$
does not satisfy the topological requirements needed to be
an $\LAR$-annulus of $\Sigma_2 \Sigma_2(P)$.
However, if we ``complete'' the faulty pattern \input{fig/E33a_nodot.tex}
by \input{fig/E33aAR_nodot.tex} in $P \cup A$ using the single white face $f$,
then $A$ becomes an $\LAR$-annulus of $P$,
as well as its image by $\Sigma_2 \Sigma_2$,
which also becomes ``completed'' by the faces in white
as shown in Equation (\ref{eq:compl}) above.
\end{exam}

\subsection{Covering properties for Arnoux-Rauzy substitutions}
\label{subsec:coverAR}
Let $\LAR$ be the set of patterns consisting of
$\input{fig/E12a.tex},
\input{fig/E12b.tex},
\input{fig/E13a.tex},
\input{fig/E13b.tex},
\input{fig/E23a.tex},
\input{fig/E23b.tex} \in \Ledge$
and of
\[
\renewcommand{\arraycolsep}{0.5mm}
\renewcommand{\arraystretch}{1.5}
\begin{array}{cclclclccclclcl}
\input{fig/E11aAR.tex} & = & [\mathbf 0, 1]^* & \cup & [(0,1,0), 1]^* & \cup & [\mathbf 0, 3]^* & \quad &
\input{fig/E11bAR.tex} & = & [\mathbf 0, 1]^* & \cup & [(0,0,1), 1]^* & \cup & [\mathbf 0, 2]^* \\
\input{fig/E22aAR.tex} & = & [\mathbf 0, 2]^* & \cup & [(1,0,0), 2]^* & \cup & [\mathbf 0, 3]^* & \quad &
\input{fig/E22bAR.tex} & = & [\mathbf 0, 2]^* & \cup & [(0,0,1), 2]^* & \cup & [\mathbf 0, 1]^* \\
\input{fig/E33aAR.tex} & = & [\mathbf 0, 3]^* & \cup & [(1,0,0), 3]^* & \cup & [\mathbf 0, 2]^* & \quad &
\input{fig/E33bAR.tex} & = & [\mathbf 0, 3]^* & \cup & [(0,1,0), 3]^* & \cup & [\mathbf 0, 1]^*.
\end{array}
\]

The goal of this section is to show that the set of patterns $\LAR$
behaves well with respect to Arnoux-Rauzy dual substitutions:
$\LAR$-covering is preserved by applying $\Sigma_1$, $\Sigma_2$ or $\Sigma_3$
(Proposition \ref{prop:covAR}),
and \emph{strong} $\LAR$-covering is also preserved
(Proposition \ref{prop:strongcovAR}),
under a pattern-avoidance condition that is always verified
in the cases that we will encounter in Section \ref{subsec:annulus} (Lemma \ref{lemm:forbidden3}).

\begin{prop}
\label{prop:covAR}
Let $P$ be an $\LAR$-covered pattern.
Then the pattern $\Sigma_i(P)$ is $\LAR$-covered for every $i \in \{1,2,3\}$.
\end{prop} 

The proof (which is detailed in \cite{BJS11}) relies on Proposition \ref{prop:coverprop}:
$\Sigma_i(P)$ is $\LAR$-covered for all $P \in \LAR$,
which can be checked by inspection.
The following stronger property will be needed in Section \ref{subsec:annulus},
to prove the annulus property for Arnoux-Rauzy substitutions.

\begin{prop}
\label{prop:strongcovAR}
Let $P$ be a strongly $\LAR$-covered pattern
that does not contain any translate of one of the three patterns
\[
\renewcommand{\arraystretch}{1.2}
\begin{array}{ccc}
\input{fig/V11a_dot.tex} & = & [(0,0,0),1]^* \cup [(0,1,1), 1]^*, \\
\input{fig/V22a_dot.tex} & = & [(0,0,0),2]^* \cup [(1,0,1), 2]^*, \\
\input{fig/V33a_dot.tex} & = & [(0,0,0),3]^* \cup [(1,1,0), 3]^*. \\
\end{array}
\]
Then $\Sigma_i(P)$ is strongly $\LAR$-covered for every $i \in \{1,2,3\}$.
\end{prop}

\begin{proof}
Let $\Sigma = \Sigma_i$ for some $i \in \{1,2,3\}$.
By Proposition \ref{prop:covAR}, $\Sigma(P)$ is $\LAR$-covered.
Now, let $X \subseteq \Sigma(P)$ be an element of $\oLedge$.
If $X$ is a translation of one of the six patterns
\input{fig/E12a_nodot.tex},
\input{fig/E12b_nodot.tex},
\input{fig/E13a_nodot.tex},
\input{fig/E13b_nodot.tex},
\input{fig/E23a_nodot.tex},
\input{fig/E23b_nodot.tex},
then the strong covering condition is trivially verified because $X$ is itself in $\oLAR$.

It remains to treat the case of the six remaining possibilities for $X$ up to translation.
We have $X \subseteq \Sigma(P)$, so it is sufficient to check that
for every pattern $X_0 \subseteq P$ such that $X \subseteq \Sigma(X_0) $, 
there exists $Y \in \oLAR$ such that $X \subseteq Y \subseteq \Sigma(X_0)$.
Moreover, since $X$ is a two-face pattern,
we can restrict to the case where $X_0$ consists of two faces only,
which leaves a finite number of cases to check for $X_0$.

Suppose that $X = \input{fig/E11a_nodot.tex}$.
The following table displays all the possibilities for $X_0$.
It has been obtained by applying Lemma \ref{lemm:preimages} to the two faces of $X$.
The first column of this table provides the index $i$ for which $\Sigma=\Sigma_i$.
\begin{center}
\renewcommand{\arraystretch}{1.2}
\begin{tabular}[h]{c|c|c}
$i$ & $X_0$ such that $X \subseteq \Sigma_i(X_0)$ & $\Sigma_i(X_0)$ \\
\hline
$2$ & \input{fig/E12a_nodot.tex}
    & \myvcenter{\begin{tikzpicture}
        [x={(-0.216506cm,-0.125000cm)}, y={(0.216506cm,-0.125000cm)}, z={(0.000000cm,0.250000cm)}]
        \definecolor{facecolor}{rgb}{0.800,0.800,0.800}
        \fill[fill=facecolor, draw=black, shift={(0,0,0)}]
        (0, 0, 0) -- (0, 0, 1) -- (1, 0, 1) -- (1, 0, 0) -- cycle;
        \fill[fill=facecolor, draw=black, shift={(0,0,0)}]
        (0, 0, 0) -- (0, 1, 0) -- (0, 1, 1) -- (0, 0, 1) -- cycle;
        \fill[fill=facecolor, draw=black, shift={(0,0,0)}]
        (0, 0, 0) -- (1, 0, 0) -- (1, 1, 0) -- (0, 1, 0) -- cycle;
        \fill[fill=facecolor, draw=black, shift={(0,1,0)}]
        (0, 0, 0) -- (0, 1, 0) -- (0, 1, 1) -- (0, 0, 1) -- cycle;
        \end{tikzpicture}}
    \\
$2$ & \input{fig/E11a_nodot.tex}
    & \input{fig/E11a_nodot.tex}
    \\
$3$ & $\myvcenter{
        \begin{tikzpicture}
        [x={(-0.216506cm,-0.125000cm)}, y={(0.216506cm,-0.125000cm)}, z={(0.000000cm,0.250000cm)}]
        \definecolor{facecolor}{rgb}{0.800,0.800,0.800}
        \fill[fill=facecolor, draw=black, shift={(0,0,0)}]
        (0, 0, 0) -- (1, 0, 0) -- (1, 1, 0) -- (0, 1, 0) -- cycle;
        \fill[fill=facecolor, draw=black, shift={(0,1,0)}]
        (0, 0, 0) -- (0, 1, 0) -- (0, 1, 1) -- (0, 0, 1) -- cycle;
        \end{tikzpicture}} = [\mathbf 0, 3]^* \cup [(0,1,0), 1]^*$
    & \myvcenter{\begin{tikzpicture}
        [x={(-0.216506cm,-0.125000cm)}, y={(0.216506cm,-0.125000cm)}, z={(0.000000cm,0.250000cm)}]
        \definecolor{facecolor}{rgb}{0.800,0.800,0.800}
        \fill[fill=facecolor, draw=black, shift={(0,0,0)}]
        (0, 0, 0) -- (0, 0, 1) -- (1, 0, 1) -- (1, 0, 0) -- cycle;
        \fill[fill=facecolor, draw=black, shift={(0,0,0)}]
        (0, 0, 0) -- (0, 1, 0) -- (0, 1, 1) -- (0, 0, 1) -- cycle;
        \fill[fill=facecolor, draw=black, shift={(0,0,0)}]
        (0, 0, 0) -- (1, 0, 0) -- (1, 1, 0) -- (0, 1, 0) -- cycle;
        \fill[fill=facecolor, draw=black, shift={(0,1,0)}]
        (0, 0, 0) -- (0, 1, 0) -- (0, 1, 1) -- (0, 0, 1) -- cycle;
        \end{tikzpicture}}
    \\
$3$ & \input{fig/V11a.tex}
    & \input{fig/E11a_nodot.tex}
\end{tabular}
\end{center}
In the first and third rows,
there is indeed a pattern $Y \in \oLAR$ such that $X \subseteq Y \subseteq \Sigma(X_0)$, namely
$Y = \input{fig/E11aAR_nodot.tex}$,
which settles the case of these two rows.
The case of the fourth row is settled directly because we have assumed that $P$ does not contain
the pattern \input{fig/V11a.tex}.
For the second row, since $P$ is itself strongly $\LAR$-covered,
$X_0$ is contained in some $Y_0\in \oLAR$
(here, \input{fig/E11aAR_nodot.tex}).
Hence, $\Sigma_2(Y_0) \subseteq \Sigma_2(P)$,
with $\Sigma_2(Y_0) \in \oLAR$, so the property holds.

The five remaining cases for $X$ can be dealt with in exactly the same way as above:
in each of the five cases, there are four preimages to check,
which are symmetrical copies of the ones in the table above.
\end{proof}

\begin{rema}
\label{rema:contreexemplestrongcovering}
The condition of avoidance of the three patterns
in Proposition \ref{prop:strongcovAR} cannot be dropped,
as illustrated by the following example:
\begin{center}
\myvcenter{%
\begin{tikzpicture}
[x={(-0.216506cm,-0.125000cm)}, y={(0.216506cm,-0.125000cm)}, z={(0.000000cm,0.250000cm)}]
\definecolor{facecolor}{rgb}{0.800,0.800,0.800}
\fill[fill=facecolor, draw=black, shift={(-1,0,2)}]
(0, 0, 0) -- (1, 0, 0) -- (1, 1, 0) -- (0, 1, 0) -- cycle;
\fill[fill=facecolor, draw=black, shift={(0,0,1)}]
(0, 0, 0) -- (0, 1, 0) -- (0, 1, 1) -- (0, 0, 1) -- cycle;
\fill[fill=facecolor, draw=black, shift={(-1,0,2)}]
(0, 0, 0) -- (0, 1, 0) -- (0, 1, 1) -- (0, 0, 1) -- cycle;
\fill[fill=facecolor, draw=black, shift={(0,0,0)}]
(0, 0, 0) -- (0, 0, 1) -- (1, 0, 1) -- (1, 0, 0) -- cycle;
\fill[fill=facecolor, draw=black, shift={(0,0,0)}]
(0, 0, 0) -- (1, 0, 0) -- (1, 1, 0) -- (0, 1, 0) -- cycle;
\fill[fill=facecolor, draw=black, shift={(-1,2,2)}]
(0, 0, 0) -- (1, 0, 0) -- (1, 1, 0) -- (0, 1, 0) -- cycle;
\fill[fill=facecolor, draw=black, shift={(-1,1,2)}]
(0, 0, 0) -- (1, 0, 0) -- (1, 1, 0) -- (0, 1, 0) -- cycle;
\fill[fill=facecolor, draw=black, shift={(0,0,0)}]
(0, 0, 0) -- (0, 1, 0) -- (0, 1, 1) -- (0, 0, 1) -- cycle;
\fill[fill=facecolor, draw=black, shift={(-1,1,2)}]
(0, 0, 0) -- (0, 1, 0) -- (0, 1, 1) -- (0, 0, 1) -- cycle;
\definecolor{facecolor}{rgb}{0.35,0.35,0.35}
\fill[fill=facecolor, draw=black, shift={(0,1,0)}]
(0, 0, 0) -- (0, 1, 0) -- (0, 1, 1) -- (0, 0, 1) -- cycle;
\fill[fill=facecolor, draw=black, shift={(0,2,1)}]
(0, 0, 0) -- (0, 1, 0) -- (0, 1, 1) -- (0, 0, 1) -- cycle;
\end{tikzpicture}%
}
\qquad $\stackrel{\Sigma_2}{\longmapsto}$ \qquad
\myvcenter{%
\begin{tikzpicture}
[x={(-0.216506cm,-0.125000cm)}, y={(0.216506cm,-0.125000cm)}, z={(0.000000cm,0.250000cm)}]
\definecolor{facecolor}{rgb}{0.800,0.800,0.800}
\fill[fill=facecolor, draw=black, shift={(-1,0,2)}]
(0, 0, 0) -- (1, 0, 0) -- (1, 1, 0) -- (0, 1, 0) -- cycle;
\fill[fill=facecolor, draw=black, shift={(0,0,1)}]
(0, 0, 0) -- (0, 1, 0) -- (0, 1, 1) -- (0, 0, 1) -- cycle;
\fill[fill=facecolor, draw=black, shift={(0,1,0)}]
(0, 0, 0) -- (0, 1, 0) -- (0, 1, 1) -- (0, 0, 1) -- cycle;
\fill[fill=facecolor, draw=black, shift={(-1,0,2)}]
(0, 0, 0) -- (0, 1, 0) -- (0, 1, 1) -- (0, 0, 1) -- cycle;
\fill[fill=facecolor, draw=black, shift={(0,1,0)}]
(0, 0, 0) -- (1, 0, 0) -- (1, 1, 0) -- (0, 1, 0) -- cycle;
\fill[fill=facecolor, draw=black, shift={(0,0,0)}]
(0, 0, 0) -- (0, 0, 1) -- (1, 0, 1) -- (1, 0, 0) -- cycle;
\fill[fill=facecolor, draw=black, shift={(0,0,0)}]
(0, 0, 0) -- (1, 0, 0) -- (1, 1, 0) -- (0, 1, 0) -- cycle;
\fill[fill=facecolor, draw=black, shift={(-1,2,2)}]
(0, 0, 0) -- (1, 0, 0) -- (1, 1, 0) -- (0, 1, 0) -- cycle;
\fill[fill=facecolor, draw=black, shift={(-1,1,2)}]
(0, 0, 0) -- (1, 0, 0) -- (1, 1, 0) -- (0, 1, 0) -- cycle;
\fill[fill=facecolor, draw=black, shift={(0,0,0)}]
(0, 0, 0) -- (0, 1, 0) -- (0, 1, 1) -- (0, 0, 1) -- cycle;
\fill[fill=facecolor, draw=black, shift={(-1,1,2)}]
(0, 0, 0) -- (0, 1, 0) -- (0, 1, 1) -- (0, 0, 1) -- cycle;
\definecolor{facecolor}{rgb}{0.35,0.35,0.35}
\fill[fill=facecolor, draw=black, shift={(0,2,0)}]
(0, 0, 0) -- (0, 1, 0) -- (0, 1, 1) -- (0, 0, 1) -- cycle;
\fill[fill=facecolor, draw=black, shift={(0,2,1)}]
(0, 0, 0) -- (0, 1, 0) -- (0, 1, 1) -- (0, 0, 1) -- cycle;
\end{tikzpicture}%
}\,.
\end{center}
Indeed, the left pattern $P$ is strongly $\LAR$-covered
and contains
\myvcenter{%
\begin{tikzpicture}
[x={(-0.216506cm,-0.125000cm)}, y={(0.216506cm,-0.125000cm)}, z={(0.000000cm,0.250000cm)}]
\definecolor{facecolor}{rgb}{0.35,0.35,0.35}
\fill[fill=facecolor, draw=black, shift={(0,0,0)}]
(0, 0, 0) -- (0, 1, 0) -- (0, 1, 1) -- (0, 0, 1) -- cycle;
\fill[fill=facecolor, draw=black, shift={(0,1,1)}]
(0, 0, 0) -- (0, 1, 0) -- (0, 1, 1) -- (0, 0, 1) -- cycle;
\end{tikzpicture}\,},
but $\Sigma_2(P)$ (depicted on the right) is not strongly $\LAR$-covered
because of the two-face edge-connected pattern shown in dark gray,
which cannot be included in any pattern $Y \in \oLAR$
such that $Y \subseteq P$.

We can observe that the pattern $P$ above has a ``hole''.
If we ``fill''
this hole with the corresponding unit face,
$P$ would no longer be strongly $\LAR$-covered.
It would be tempting to try to put some topological restrictions on $P$
to get rid of the  assumption of avoidance of the three patterns in Proposition \ref{prop:strongcovAR},
but we will actually need to apply Proposition \ref{prop:strongcovAR}
to some patterns that are homeomorphic to topological annuli (as $P$ is),
in the proof of Lemma \ref{lemm:annulus_induction}.
\end{rema}

In order to be able to apply Proposition \ref{prop:strongcovAR} to patterns
obtained by iterating $\Sigma_1$, $\Sigma_2$ and $\Sigma_3$ from $\mcU$,
we need the following lemma,
which states that none of the three above-mentioned forbidden patterns
appears in the iterations.

\begin{lemm}
\label{lemm:forbidden3}
Let $\Sigma$ be an arbitrary product of $\Sigma_1$, $\Sigma_2$ and $\Sigma_3$.
The discrete plane $\Sigma(\Gamma_{(1,1,1)})$ contains
no translate of any of the three patterns
\input{fig/V11a.tex}, \input{fig/V22a.tex}, \input{fig/V33a.tex}.
\end{lemm}

\begin{proof}
We prove the result by induction on the size of the product $\Sigma$.
First, we observe that $\Gamma_{(1,1,1)}$ does not contain any of the three patterns in question:
consider for instance a translate of the pattern
$\input{fig/V11a.tex} = [\bfx, 1]^* \cup [\bfx + (0,1,1), 1]^*$,
for $\bfx \in \bbZ^3$.
We have
$\langle \bfx + (0,1,1), (1,1,1) \rangle
    = \langle \bfx, (1,1,1) \rangle + 2
    \geq \langle \bfe_1, (1,1,1) \rangle$,
which contradicts Proposition \ref{prop:arithplane}.
The same reasoning holds for the two remaining patterns 
\input{fig/V22a.tex} and \input{fig/V33a.tex} for the initialization of the induction.

In the table below,
we have listed all the preimages of the three patterns (in light gray) obtained by applying Lemma \ref{lemm:preimages},
together with their only possible ``completion'' within a discrete plane (in dark gray).
This completion is deduced from the arithmetic description of discrete planes provided by Proposition \ref{prop:arithplane}.
\begin{center}
\renewcommand{\arraystretch}{1.7}
\begin{tabular}{c|c}
\input{fig/LVBAD1.tex}
\end{tabular}
\qquad
\begin{tabular}{c|c}
\input{fig/LVBAD2.tex}
\end{tabular}
\qquad
\begin{tabular}{c|c}
\input{fig/LVBAD3.tex}
\end{tabular}
\end{center}
Let us work by contradiction and assume that a translate of one of the three patterns
appears in $\Sigma(\Gamma_{(1,1,1)})$ but not in $\Sigma'(\Gamma_{(1,1,1)})$,
where $\Sigma=\Sigma_i  \circ  \Sigma'$, for some $i \in \{1,2,3\}$.
The above table, together with the injectivity of $\Sigma_i$ (see Proposition~\ref{prop:imgplane}),
implies that a translate of one of the three patterns must appear in $\Sigma'(\Gamma_{(1,1,1)})$,
which yields the desired contradiction.
\end{proof}

\section{Main results}
\label{sec:mainresults}
We now use the results of Section \ref{sec:cov} to prove
the annulus property for Arnoux-Rauzy substitutions,
from which we will deduce our main theorem and some other dynamical consequences.

\subsection{The annulus property for Arnoux-Rauzy substitutions}
\label{subsec:annulus}

The aim of this section is to prove the following theorem.
The proof (by induction) relies on Lemma \ref{lemm:annulus_base} (base case)
and on Lemma \ref{lemm:annulus_induction} (induction step).

\begin{theo}
\label{theo:annulus}
Let $\Sigma$ be a product of $\Sigma_1$, $\Sigma_2$ and $\Sigma_3$
in which each substitution appears at least once.
Then, $\Sigma^{n+2}(\mcU) \setminus \Sigma^n(\mcU)$ is an $\LAR$-annulus of $\Sigma^n(\mcU)$
in the discrete plane $\Sigma^{n+2}(\Gamma_{(1,1,1)})$, for all $n \in \bbN$.
\end{theo}

\begin{proof}
We prove the result by induction on $n$.
The case $n=0$ is settled by Lemma \ref{lemm:annulus_base}.
Now, assume that the induction property holds for some $n \in \bbN$.
According to Proposition \ref{prop:imgplane},
$\Sigma^{n+2}(\mcU)$ is contained in the discrete plane $\Sigma^{n+2}(\Gamma_{(1,1,1)})$.
Hence, Lemma \ref{lemm:forbidden3} allows us to apply Lemma \ref{lemm:annulus_induction}
to deduce that $\Sigma^{n+3}(\mcU) \setminus \Sigma^{n+1}(\mcU)$
is an $\LAR$-annulus of $\Sigma^{n+1}(\mcU)$.
\end{proof}

\begin{lemm}
\label{lemm:annulus_base}
Let $\Sigma$ be a product of $\Sigma_1$, $\Sigma_2$ and $\Sigma_3$
in which each substitution appears at least once.
Then, $\Sigma^2(\mcU) \setminus \mcU$ is an $\LAR$-annulus of $\mcU$.
\end{lemm}

\begin{figure}[ht]
\centering
\begin{tikzpicture}[x={(-0.866cm,-0.5cm)}, y={(0.866cm,-0.5cm)}, z={(0cm,1cm)}]
%[x=1cm,y=0.85cm]
%\node at (0,-0.5) (0)   {\includegraphics[scale=0.85]{fig/g0.pdf}};

%\node at (-4,-2) (1)    {\includegraphics[scale=0.85]{fig/g1.pdf}};
%\node at (0,-2)  (2)    {\includegraphics[scale=0.85]{fig/g2.pdf}};
%\node at (4,-2)  (3)    {\includegraphics[scale=0.85]{fig/g3.pdf}};

%\node at (-3,-4) (12)   {\includegraphics[scale=0.85]{fig/g12.pdf}};
%\node at (-5,-4) (13)   {\includegraphics[scale=0.85]{fig/g13.pdf}};
%\node at (-1,-4) (21)   {\includegraphics[scale=0.85]{fig/g21.pdf}};
%\node at (1,-4)  (23)   {\includegraphics[scale=0.85]{fig/g23.pdf}};
%\node at (5,-4)  (31)   {\includegraphics[scale=0.85]{fig/g31.pdf}};
%\node at (3,-4)  (32)   {\includegraphics[scale=0.85]{fig/g32.pdf}};

%\node at (-3,-6) (123)  {\includegraphics[scale=0.85]{fig/g123.pdf}};
%\node at (-5,-6) (132)  {\includegraphics[scale=0.85]{fig/g132.pdf}};
%\node at (-1,-6) (213)  {\includegraphics[scale=0.85]{fig/g213.pdf}};
%\node at (1,-6)  (231)  {\includegraphics[scale=0.85]{fig/g231.pdf}};
%\node at (5,-6)  (312)  {\includegraphics[scale=0.85]{fig/g312.pdf}};
%\node at (3,-6)  (321)  {\includegraphics[scale=0.85]{fig/g321.pdf}};

%\node at (-3,-8) (1231) {\includegraphics[scale=0.85]{fig/g1231.pdf}};
%\node at (-5,-8) (1321) {\includegraphics[scale=0.85]{fig/g1321.pdf}};
%\node at (-1,-8) (2132) {\includegraphics[scale=0.85]{fig/g2132.pdf}};
%\node at (1,-8)  (2312) {\includegraphics[scale=0.85]{fig/g2312.pdf}};
%\node at (5,-8)  (3123) {\includegraphics[scale=0.85]{fig/g3123.pdf}};
%\node at (3,-8)  (3213) {\includegraphics[scale=0.85]{fig/g3213.pdf}};

\node at (0,0,0) (0)     {\includegraphics[scale=0.8]{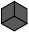}};

\node at (2,0,0) (1)     {\includegraphics[scale=0.8]{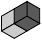}};
\node at (0,2,0) (2)     {\includegraphics[scale=0.8]{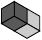}};
\node at (0,0,2) (3)     {\includegraphics[scale=0.8]{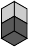}};

\node at (2,0,-2) (12)   {\includegraphics[scale=0.8]{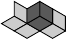}};
\node at (2,-2,0) (13)   {\includegraphics[scale=0.8]{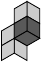}};
\node at (0,2,-2) (21)   {\includegraphics[scale=0.8]{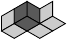}};
\node at (-2,2,0) (23)   {\includegraphics[scale=0.8]{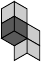}};
\node at (0,-2,2) (31)   {\includegraphics[scale=0.8]{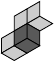}};
\node at (-2,0,2) (32)   {\includegraphics[scale=0.8]{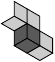}};

\node at (4,0,-2) (123)  {\includegraphics[scale=0.8]{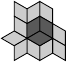}};
\node at (4,-2,0) (132)  {\includegraphics[scale=0.8]{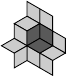}};
\node at (0,4,-2) (213)  {\includegraphics[scale=0.8]{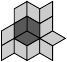}};
\node at (-2,4,0) (231)  {\includegraphics[scale=0.8]{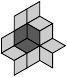}};
\node at (0,-2,4) (312)  {\includegraphics[scale=0.8]{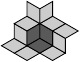}};
\node at (-2,0,4) (321)  {\includegraphics[scale=0.8]{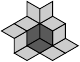}};

\node at (6,0,-2) (1231) {\includegraphics[scale=0.8]{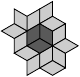}};
\node at (6,-2,0) (1321) {\includegraphics[scale=0.8]{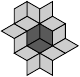}};
\node at (0,6,-2) (2132) {\includegraphics[scale=0.8]{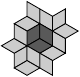}};
\node at (-2,6,0) (2312) {\includegraphics[scale=0.8]{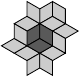}};
\node at (0,-2,6) (3123) {\includegraphics[scale=0.8]{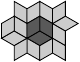}};
\node at (-2,0,6) (3213) {\includegraphics[scale=0.8]{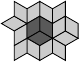}};

\path[->, thick] (0) edge node [below] {$1$} (1);
\path[->, thick] (0) edge node [above] {$2$} (2);
\path[->, thick] (0) edge node [left] {$3$} (3);

\path[->, thick] (1) edge node [right] {$2$} (12);
\path[->, thick] (1) edge node [above] {$3$} (13);
\path[->, thick] (2) edge node [left] {$1$} (21);
\path[->, thick] (2) edge node [above] {$3$} (23);
\path[->, thick] (3) edge node [below] {$1$} (31);
\path[->, thick] (3) edge node [below] {$2$} (32);

\path[->, thick] (12) edge [bend left] node [above] {$1$} (21);
\path[->, thick] (13) edge [bend left] node [left] {$1$} (31);
\path[->, thick] (21) edge [bend left] node [below] {$2$} (12);
\path[->, thick] (23) edge [bend left] node [left] {$2$} (32);
\path[->, thick] (31) edge [bend left] node [right] {$3$} (13);
\path[->, thick] (32) edge [bend left] node [right] {$3$} (23);

\path[->, thick] (12) edge node [below] {$3$} (123);
\path[->, thick] (13) edge node [above] {$2$} (132);
\path[->, thick] (21) edge node [below] {$3$} (213);
\path[->, thick] (23) edge node [above] {$1$} (231);
\path[->, thick] (31) edge node [left] {$2$} (312);
\path[->, thick] (32) edge node [right] {$1$} (321);

\path[->, thick] (123) edge [bend left] node [left] {$2$} (132);
\path[->, thick] (132) edge [bend left] node [right] {$3$} (123);
\path[->, thick] (213) edge [bend left] node [left] {$1$} (231);
\path[->, thick] (231) edge [bend left] node [right] {$3$} (213);
\path[->, thick] (312) edge [bend left] node [above] {$1$} (321);
\path[->, thick] (321) edge [bend left] node [below] {$2$} (312);

\path[->, thick] (123) edge node [below] {$1$} (1231);
\path[->, thick] (132) edge node [above] {$1$} (1321);
\path[->, thick] (213) edge node [below] {$2$} (2132);
\path[->, thick] (231) edge node [above] {$2$} (2312);
\path[->, thick] (312) edge node [left] {$3$} (3123);
\path[->, thick] (321) edge node [right] {$3$} (3213);

% loops
%\draw[->, thick] (1) .. controls +(110:10mm) and +(150:10mm) .. node [above left] {$1$} (1);
%\draw[->, thick] (2) .. controls +(60:10mm) and +(20:10mm) .. node [above right] {$2$} (2);
%\draw[->, thick] (3) .. controls +(60:10mm) and +(20:10mm) .. node [above right] {$3$} (3);
%
%\draw[->, thick] (12) .. controls +(160:10mm) and +(200:10mm) .. node [left] {$2$} (12);
%\draw[->, thick] (13) .. controls +(160:10mm) and +(200:10mm) .. node [left] {$3$} (13);
%\draw[->, thick] (21) .. controls +(160:10mm) and +(200:10mm) .. node [left] {$1$} (21);
%\draw[->, thick] (23) .. controls +(-20:10mm) and +(20:10mm) .. node [right] {$3$} (23);
%\draw[->, thick] (31) .. controls +(-20:10mm) and +(20:10mm) .. node [right] {$1$} (31);
%\draw[->, thick] (32) .. controls +(-20:10mm) and +(20:10mm) .. node [right] {$2$} (32);

%dges
%\path[->, thick] (1231) edge [bend right]  node [above] {$2$} (2312);
%\path[->, thick] (1231) edge [bend right]  node [above] {$3$} (3213);
%\path[->, thick] (1321) edge [bend right]  node [above] {$2$} (2312);
%\path[->, thick] (1321) edge [bend right]  node [above] {$3$} (3213);
\end{tikzpicture}
\caption{Each vertex has three outgoing edges with distinct labels $1$, $2$ or $3$.
(The $9$ loops and the $18$ outgoing edges of the extremal vertices are not drawn.)}
\label{fig:ARgraph}
\end{figure}

\begin{proof}
Let $w = i_1 \cdots i_n \in \{1,2,3\}^n$ such that $\Sigma = \Sigma_{i_1} \cdots \Sigma_{i_n}$.
This proof is a case study, which we formalize as the study
of the directed graph in Fig. \ref{fig:ARgraph}.
The vertices of this finite  graph are patterns.
This graph is such that for every edge
$P \stackrel{i}{\longrightarrow} Q$ we have $Q \subseteq \Sigma_i(P)$,
and every vertex has three outgoing edges with distinct labels $1$, $2$ and $3$
(for the sake of clarity, loops are not drawn in Fig. \ref{fig:ARgraph}).

It can be checked that one of the six extremal vertices is always reached when
following the path of labelled edges given by $w^2=ww$ starting at the top vertex $\mcU$.

The result follows because each extremal vertex
is a pattern that consists of $\mcU$ and an $\LAR$-annulus of $\mcU$,
and because each outgoing edge of an extremal vertex points to an extremal vertex
(these edges are not drawn in Fig. \ref{fig:ARgraph}),
so applying $\Sigma$ twice guarantees us to reach an $\LAR$-annulus of $\mcU$.
\end{proof}

Let us note that the existence of a finite graph such as the one in Fig. \ref{fig:ARgraph}
that provides the description of the possible annuli surrounding the unit cube
is the key ingredient of the proof.
There is no reason \emph{a priori} for such a finite graph to exist for a given set of substitutions.

\begin{lemm}
\label{lemm:annulus_induction}
Let $\Gamma$ be a discrete plane that does not contain any of the patterns
\input{fig/V11a.tex}, \input{fig/V22a.tex}, \input{fig/V33a.tex}.
Let $A \subseteq \Gamma$ be an $\LAR$-annulus of a pattern $P \subseteq \Gamma$,
and let $\Sigma = \Sigma_i$ for some $i \in \{1,2,3\}$.
Then $\Sigma(A)$ is an $\LAR$-annulus of $\Sigma(P)$.
\end{lemm}

\begin{proof}
We must prove the following:
\begin{enumerate}
  \item $\Sigma(P)$, $\Sigma(A) \cup \Sigma(P)$ and $\Gamma \setminus (\Sigma(A) \cup \Sigma(P))$ are $\LAR$-covered;
  \item $\Sigma(A)$ is strongly $\LAR$-covered;
  \item $\Sigma(A)$ and $\Sigma(P)$ have no face in common;
  \item $\Sigma(P) \cap \overline{\Sigma(\Gamma) \setminus (\Sigma(P) \cup \Sigma(A))} = \varnothing$.
\end{enumerate}
Conditions (\ref{defi:anneauprop1}) and (\ref{defi:anneauprop3}) hold
by Propositions \ref{prop:coverprop} and \ref{prop:imgplane} respectively,
and (\ref{defi:anneauprop2}) holds thanks to Proposition \ref{prop:strongcovAR}.
It remains to prove that (\ref{defi:anneauprop4}) also holds.

We will work by contradiction and we suppose that (\ref{defi:anneauprop4}) does not hold.
This implies that there exist faces
\[
f \in P, \qquad
g \in \Gamma \setminus (A \cup P), \qquad
f' \in \Sigma(f), \qquad
g' \in \Sigma(g)
\]
such that $f'$ and $g'$ have a nonempty intersection.
Note that $f \cup g$ is disconnected because $P$ and $\overline{\Gamma \setminus (P \cup A)}$
have empty intersection by hypothesis.

\begin{figure}[ht]
\centering
    \myvcenter{\includegraphics[scale=0.85]{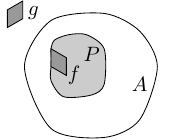}} \qquad
    \myvcenter{\includegraphics[scale=0.85]{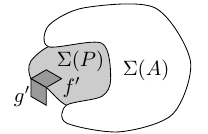}} \qquad
    \myvcenter{\includegraphics[scale=0.85]{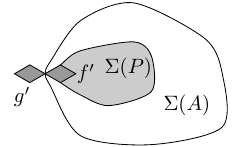}}
\caption[]{The patterns $A \cup P$ (left) and $\Sigma(A) \cup \Sigma(P)$ (middle or right),
    under the assumption that Condition (\ref{defi:anneauprop4}) does not hold.}
\label{fig:lemcases}
\end{figure}

The strategy of the proof is as follows:
we check all the possible patterns $f \cup g$ and $f' \cup g'$ as above,
and for each case we derive a contradiction.
This can be done by inspection of a finite number of cases.
Indeed, there are $36$ possibilities for $f' \cup g'$ up to translation
(the number of connected two-face patterns that share a vertex or an edge),
and each of these patterns has a finite number of two-face preimages
given by Lemma \ref{lemm:preimages}.

The first such case is the following:
\[
\left\{
\begin{array}{ccccl}
f \cup g
    & = & \myvcenter{%
    \begin{tikzpicture}
    [x={(-0.216506cm,-0.125000cm)}, y={(0.216506cm,-0.125000cm)}, z={(0.000000cm,0.250000cm)}]
    \definecolor{facecolor}{rgb}{0.8,0.8,0.8}
    \fill[fill=facecolor, draw=black, shift={(0,0,0)}]
    (0, 0, 0) -- (0, 1, 0) -- (0, 1, 1) -- (0, 0, 1) -- cycle;
    \fill[fill=facecolor, draw=black, shift={(1,-1,1)}]
    (0, 0, 0) -- (0, 1, 0) -- (0, 1, 1) -- (0, 0, 1) -- cycle;
    \end{tikzpicture}}
    & = & [\mathbf 0, 1]^* \cup [(1, -1, 1), 1]^* \\
f' \cup g'
    & = & \myvcenter{%
    \begin{tikzpicture}
    [x={(-0.216506cm,-0.125000cm)}, y={(0.216506cm,-0.125000cm)}, z={(0.000000cm,0.250000cm)}]
    \definecolor{facecolor}{rgb}{0.8,0.8,0.8}
    \fill[fill=facecolor, draw=black, shift={(0,0,0)}]
    (0, 0, 0) -- (0, 1, 0) -- (0, 1, 1) -- (0, 0, 1) -- cycle;
    \fill[fill=facecolor, draw=black, shift={(-1,1,-1)}]
    (0, 0, 0) -- (0, 0, 1) -- (1, 0, 1) -- (1, 0, 0) -- cycle;
    \end{tikzpicture}}
    & = & [\mathbf 0, 2]^* \cup [(1, -1, 1), 1]^*.
\end{array}
\right.
\]
%    \ \subseteq \
%    \Sigma_1(f \cup g) \ = \ \myvcenter{%
%    \begin{tikzpicture}
%    [x={(-0.216506cm,-0.125000cm)}, y={(0.216506cm,-0.125000cm)}, z={(0.000000cm,0.250000cm)}]
%    \definecolor{facecolor}{rgb}{0.8,0.8,0.8}
%    \fill[fill=facecolor, draw=black, shift={(0,0,0)}]
%    (0, 0, 0) -- (0, 1, 0) -- (0, 1, 1) -- (0, 0, 1) -- cycle;
%    \fill[fill=facecolor, draw=black, shift={(0,0,0)}]
%    (0, 0, 0) -- (0, 0, 1) -- (1, 0, 1) -- (1, 0, 0) -- cycle;
%    \fill[fill=facecolor, draw=black, shift={(0,0,0)}]
%    (0, 0, 0) -- (1, 0, 0) -- (1, 1, 0) -- (0, 1, 0) -- cycle;
%    \fill[fill=facecolor, draw=black, shift={(1,-1,1)}]
%    (0, 0, 0) -- (0, 1, 0) -- (0, 1, 1) -- (0, 0, 1) -- cycle;
%    \fill[fill=facecolor, draw=black, shift={(1,-1,1)}]
%    (0, 0, 0) -- (0, 0, 1) -- (1, 0, 1) -- (1, 0, 0) -- cycle;
%    \fill[fill=facecolor, draw=black, shift={(1,-1,1)}]
%    (0, 0, 0) -- (1, 0, 0) -- (1, 1, 0) -- (0, 1, 0) -- cycle;
%    \end{tikzpicture}}
By Proposition \ref{prop:arithplane},
the pattern $f \cup g$ can only be ``completed'' in the following way
within the discrete plane $\Gamma$:
\[
f \cup g \cup X \ = \ \myvcenter{%
\begin{tikzpicture}
[x={(-0.216506cm,-0.125000cm)}, y={(0.216506cm,-0.125000cm)}, z={(0.000000cm,0.250000cm)}]
\definecolor{facecolor}{rgb}{0.8,0.8,0.8}
\fill[fill=facecolor, draw=black, shift={(0,0,0)}]
(0, 0, 0) -- (0, 1, 0) -- (0, 1, 1) -- (0, 0, 1) -- cycle;
\fill[fill=facecolor, draw=black, shift={(1,-1,1)}]
(0, 0, 0) -- (0, 1, 0) -- (0, 1, 1) -- (0, 0, 1) -- cycle;
\definecolor{facecolor}{rgb}{0.35,0.35,0.35}
\fill[fill=facecolor, draw=black, shift={(0,0,0)}]
(0, 0, 0) -- (0, 0, 1) -- (1, 0, 1) -- (1, 0, 0) -- cycle;
\fill[fill=facecolor, draw=black, shift={(0,0,1)}]
(0, 0, 0) -- (0, 0, 1) -- (1, 0, 1) -- (1, 0, 0) -- cycle;
\end{tikzpicture}%
}\,,
\]
where $X = [\mathbf 0, 2]^* \cup [(0, 0, 1), 2]^*$ is shown in dark gray.
Moreover, we have $X \subseteq A$ because Condition (\ref{defi:anneauprop4}) for $A$ and $P$ would fail otherwise.
However, this situation is impossible:
$X \in \oLedge$ but there cannot exist a pattern $Y \in \oLAR$ such that $X \subseteq Y \subseteq A$
because $Y$ would overlap with $f$ or $g$, which both are not in $A$.
This implies that $A$ is not strongly $\LAR$-covered, which is a contradiction.

The same reasoning applies to all the other possible cases for $f'$ and $g'$,
as shown by the table of Fig. \ref{fig:tablescases}.
(By assumption, we have not considered the (problematic) patterns
\input{fig/V11a.tex}, \input{fig/V22a.tex}, \input{fig/V33a.tex}.)
\end{proof}

\begin{figure}[ht]
\centering
\renewcommand{\arraystretch}{1.65}
\input{fig/badcasesAR1.tex}
\qquad
\input{fig/badcasesAR2.tex}
\qquad
\input{fig/badcasesAR3.tex}
\caption{Disconnected preimages of connected two-face patterns by
  $\Sigma_1$ (left), $\Sigma_2$ (middle) and $\Sigma_3$ (right).
  The only possible completion of $f \cup g$ is in dark gray,
  and $f \cup g$, $f' \cup g'$ are in light gray.}
\label{fig:tablescases}
\end{figure}

\subsection{Dynamical consequences}
\label{subsec:dynconseq}

Our main theorems now follow easily from Theorem \ref{theo:annulus}.

\begin{theo}
\label{theo:growingballs}
Let $\sigma$ be a finite product of the Arnoux-Rauzy substitutions $\sigma_1$, $\sigma_2$, $\sigma_3$
in which each substitution appears at least once, and let $\Sigma = \EOSS$.
Then, for all $R \geq 0$, there exists $n \geq 0$ such that the ball of radius $R$ centered at the origin
is contained in $\pic(\Sigma^n(\mcU))$.
\end{theo}

\begin{proof}
Let $n \geq 0$ and $\Gvb$ be the discrete plane associated with the contracting plane of $\sigma$
(we have $\Sigma^n(\mcU) \subseteq \Gvb$ for all $n \geq 0$).
By Theorem \ref{theo:annulus},
the pattern $\Sigma^{2n}(\mcU)$ consists of $\mcU$ and $n$ disjoint $\LAR$-annuli $A_1, \ldots, A_n$
such that $A_1$ is an $\LAR$-annulus of $\mcU$ and $A_k$ is an $\LAR$-annulus of
$\mcU \cup A_1 \cup \cdots A_{k-1}$ for all $2 \leq k \leq n$,
in the discrete plane $\Gvb$.

We define the \emph{combinatorial radius of} $\Sigma^{2n}(\mcU)$
as the smallest edge-connected discrete path of unit faces
from a face of $\mcU$ to a face of $\Gvb \setminus \Sigma^{2n}(\mcU)$.
Condition (\ref{defi:anneauprop1}) of Definition \ref{defi:annulus}
ensures that $\Sigma^{2n}(\mcU)$ is simply connected
(because its complement is $\LAR$-covered and the patterns of $\LAR$ are edge-connected),
so Conditions (\ref{defi:anneauprop3}) and (\ref{defi:anneauprop4}) ensure that each annulus
adds at least $1$ to the combinatorial radius of $\Sigma^{2n}(\mcU)$
(again, because the patterns of $\LAR$ are edge-connected).

It follows that the combinatorial radius of $\Sigma^{2n}(\mcU)$ is at least $n$,
which implies that $\pic(\Sigma^n(\mcU))$ contains arbitrarily large balls when $n$ goes to infinity.
\end{proof}

We can now apply Theorem \ref{theo:puredis} to obtain Theorem \ref{theo:intro}:

\begin{coro}
Let $\sigma$ be a finite product of the Arnoux-Rauzy substitutions $\sigma_1$, $\sigma_2$ and $\sigma_3$
in which each substitution appears at least once.
Then, the symbolic dynamical system generated by $\sigma$
is measurably conjugate to a toral translation,
\emph{i.e.}, it has pure discrete spectrum.
\end{coro}

We also derive further topological properties of the associated Rauzy fractal.

\begin{theo}
\label{theo:RFprop}
Let $\sigma$ be a finite product of the Arnoux-Rauzy substitutions $\sigma_1$, $\sigma_2$ and $\sigma_3$
in which each substitution appears at least once.
The Rauzy fractal associated with
$\sigma$ is connected and $\mathbf 0$ is an inner point.
\end{theo}

The connectedness has been established in \cite{BJS11},
and the fact that $\mathbf 0$ is an inner point
comes from  Theorem \ref{theo:growingballs} together with \cite{ST10}.

More generally, we have the following ``$S$-adic version'' of Theorem \ref{theo:growingballs}.

\begin{theo}
\label{theo:growingballsSadic}
Let $u$ be an Arnoux-Rauzy sequence defined as 
\[
u = \lim_{n \rightarrow \infty}   \sigma_{i_0} \circ \cdots \circ \sigma _{i_n},
\]
where $\sigma_i$ ($i=1,2,3$) are the Arnoux-Rauzy substitutions,
and the sequence $(i_n)_{n \in \bbN} \in \{1,2,3\}^\bbN$ takes each value in $\{1,2,3\}$ an infinite number of times.
Then, for all $R \geq 0$, there exists $n \geq 0$ such that the ball of radius $R$ centered at the origin
is contained in
\[
\pic(\EOS(\sigma_{i_0}) \circ \cdots \circ \EOS(\sigma_{i_n})(\mcU)).
\]
\end{theo}

\begin{proof}
The proof follows the same lines as the proof of Theorem~\ref{theo:growingballs}.
Let $(n_k)_{k \in \bbN}$ be an increasing sequence of integers such that the word
$i_{n _k +1} \cdots i_{n_{k+1}}$ contains each of the letters $1,2,3$ for every $k$.
The role played by the substitution $\sigma$
will be played here by $\sigma_{i_{{n_k}+1}} \circ \cdots \circ \sigma_{i_{n_{k+1}}}$,
and the role played by $\Sigma$ 
will be played  by
$\EOS(\sigma_{i_{n_{k+1}}}) \circ \cdots \circ \EOS (\sigma_{i_{n_k+1}})
 = \EOS(\sigma_{i_{{n_k}+1}} \circ \cdots \circ \sigma_{i_{n_{k+1}}})$.
Let
\[
\Sigma_k = \EOS  (\sigma_{i_0})  \circ  \cdots  \circ   \EOS  (\sigma_{i_{n_k}}).
\]
One proves exactly as for Theorem \ref{theo:annulus}
that $\Sigma_{k+2}(\mcU) \setminus \Sigma_k(\mcU)$ is an $\LAR$-annulus of $\Sigma_k(\mcU)$
in the discrete plane $\Sigma_{k+2}(\Gamma_{(1,1,1)})$, for all $k \in \bbN$.
Indeed, by definition of the sequence $(n_k)_{k \in \bbN}$, 
in the proof of Lemma \ref{lemm:annulus_base},
any path labelled by $i_{n_k+1} \cdots i_{n_{k+2}}$ starting at the top vertex $\mcU$ reaches 
one of the six  bottom vertices.
\end{proof}

Such a result can be seen as a generation method
(in the vein of \cite{FerBrun, BLPP11}) for a discrete plane
whose normal vector is in some sense dual to the vector of
letter frequencies associated with the infinite word $u$.
It can also be considered as a first step towards
the definition of Rauzy fractals and Rauzy tilings in the $S$-adic framework. 

\appendix
\section{Proof of Lemma \ref{lemm:preimages}}
\label{app:preimages}
\begin{proof}
We treat the case of $\Sigma_1$ only,
the other cases being analogous.
The images of the faces by $\Sigma_1$ are given by
\[
\begin{array}{rcl}
\big[\svect{x}{y}{z}, 1\big]^*
  & \mapsto &
  \svect{x-y-z}{y}{z} + \big[\svect{0}{0}{0}, 1\big]^* \cup \big[\svect{0}{0}{0}, 2\big]^* \cup \big[\svect{0}{0}{0}, 3\big]^* \\
\big[\svect{x}{y}{z}, 2\big]^*
  & \mapsto & \svect{x-y-z}{y}{z} + \big[\svect{1}{0}{0}, 2\big]^* \\
\big[\svect{x}{y}{z}, 3\big]^*
  & \mapsto & \svect{x-y-z}{y}{z} + \big[\svect{1}{0}{0}, 3\big]^*
\end{array}.
\]
A face of type $1$ appears only once in the image of a face of type $1$.
If $[(x,y,z), 1]^*$ belongs to $\Sigma_1([(x_0, y_0, z_0), 1]^*)$,
then
\[
  \left \{
    \begin{array}{rcl}
    x & = & x_0 - y_0 - z_0 \\
    y & = & y_0 \\
    z & = & z_0
    \end{array}
  \right.
\quad \Longrightarrow \quad
  \left \{
    \begin{array}{rcl}
    x_0 & = & x + y + z \\
    y_0 & = & y \\
    z_0 & = & z
    \end{array}
  \right..
\]
Hence
\[
\Sigma_1^{-1}([(x,y,z), 1]^*) \ = \
  \left[ \svect{x+y+z}{y}{z}, 1 \right]^*.
\]
A face of type $2$ appears once in the image of a face of type $1$
and once in the image of a face of type $2$.
If $[(x,y,z), 2]^* \subseteq \Sigma_1([(x_0, y_0, z_0), 1]^*)$,
then
\[
  \left \{
    \begin{array}{rcl}
    x & = & x_0 - y_0 - z_0 \\
    y & = & y_0 \\
    z & = & z_0
    \end{array}
  \right.
\quad \Longrightarrow \quad
  \left \{
    \begin{array}{rcl}
    x_0 & = & x + y + z \\
    y_0 & = & y \\
    z_0 & = & z
    \end{array}
  \right..
\]
If $[(x,y,z), 2]^* \subseteq \Sigma_1([(x_0, y_0, z_0), 2]^*)$,
then
\[
  \left \{
    \begin{array}{rcl}
    x & = & x_0 - y_0 - z_0 + 1\\
    y & = & y_0 \\
    z & = & z_0
    \end{array}
  \right.
\quad \Longrightarrow \quad
  \left \{
    \begin{array}{rcl}
    x_0 & = & x + y + z - 1\\
    y_0 & = & y \\
    z_0 & = & z
    \end{array}
  \right..
\]
Hence
\[
\Sigma_1^{-1}([(x,y,z), 2]^*) \ = \
  \left[ \svect{x+y+z}{y}{z}, 1 \right]^*
  \cup
  \left[ \svect{x+y+z-1}{y}{z}, 2 \right]^*.
\]
A face of type $3$ appears once in the image of a face of type $1$
and once in the image of a face of type $3$.
In the same way as for a face of type $2$, we get
\[
\Sigma_1^{-1}([(x,y,z), 3]^*) \ = \
  \left[ \svect{x+y+z}{y}{z}, 1 \right]^*
  \cup
  \left[ \svect{x+y+z-1}{y}{z}, 3 \right]^*.
\]
\end{proof}
\begin{ackn}
Many of the computations found in this article were carried out using the Sage mathematics software system \cite{sage}.
This work was supported by the \emph{Agence Nationale de la Recherche} through contract ANR-2010-BLAN-0205-01.
\end{ackn}

\bibliographystyle{amsalpha}
\bibliography{biblio}
\end{document}